\documentclass[11pt]{article} 
\usepackage[english]{babel}

\newtheorem{definition}{Definition}
\newtheorem{theorem}{Theorem}

\usepackage{amssymb}
\usepackage{amsmath}
\def\ds{\displaystyle}
\def\ph{\phi}
\def\dt{{\delta t}}
\def\dx{{\delta x}}
\def\dy{{\delta y}}

\def\lx{{\lambda_x}}
\def\ly{{\lambda_y}}
\def\eps{\varepsilon}
\def\rmr{{\mathrm{r}}}
\def\rms{{\mathrm{s}}}
\def\bfB{{\mathbf{B}}}
\def\bfD{{\mathbf{D}}}
\def\bfE{{\mathbf{E}}}
\def\bfJ{{\mathbf{J}}}
\def\bfP{{\mathbf{P}}}
\def\bfj{{\mathbf{j}}}
\def\bfx{{\mathbf{x}}}
\def\bfxi{{\boldsymbol{\xi}}}
\def\bbC{{\mathbb{C}}}

\def\bbR{{\mathbb{R}}}
\def\calB{{\mathcal{B}}}
\def\calD{{\mathcal{D}}}
\def\calE{{\mathcal{E}}}
\def\calJ{{\mathcal{J}}}
\def\calP{{\mathcal{P}}}
\def\Id{{\mathrm{Id}}}
\def\curl{{\operatorname{\ curl\ }}}
\def\transp{^t}
\begin{document}

\title{Von Neumann Stability Analysis of Finite Difference Schemes \\
for Maxwell--Debye and Maxwell--Lorentz Equations}

\author{Brigitte Bid\'egaray-Fesquet \\
\small Laboratoire Jean-Kuntzmann, Grenoble University and UMR CNRS 5224\\
\small B.P. 53, 38041 Grenoble Cedex 9, France}

\date{}
\maketitle

\begin{abstract}
This technical report yields detailed calculations of the paper
\cite{BidegarayFesquet05c} which have been however automated since (see
\cite{BidegarayFesquet06c}). It deals with the stability analysis of various
finite difference schemes for Maxwell--Debye and Maxwell--Lorentz
equations. This work gives a systematic and rigorous continuation to
Petropoulos previous work \cite{Petropoulos94a}.

\end{abstract}

\section{Introduction}

We address the stability study of finite difference schemes for Maxwell--Debye
and Maxwell--Lorentz models. To this aim we selected the same schemes as those
already studied by Petropoulos \cite{Petropoulos94a}, who after having 
correctly defined characteristic polynomials associated to each scheme, merely 
computed its roots with a numerical algorithm. This implies having to specify
values for the physical parameters which occur in the models as well for the
time and space steps chosen for the discretization. The analysis has therefore
to be carried out anew for each new material or discretization. We perform here
a von Neumann analysis on the characteristic polynomials in their literal form,
which yields once and for all stability conditions which are valid for all
materials.

\subsection{Maxwell--Debye and Maxwell--Lorentz Models}

Le us consider Maxwell equations without magnetisation
\begin{equation}
\label{Maxwell}
\begin{array}{lrcl}
\textrm{(Faraday)} & 
\ds\partial_t \bfB(t,\bfx) & = & \ds- \curl \bfE(t,\bfx), \\
\textrm{(Amp\`ere)} &
\ds\partial_t \bfD(t,\bfx) & = & \ds\frac1{\mu_0} \curl \bfB(t,\bfx),
\end{array}
\end{equation}
where $\bfx\in\bbR^N$. This system is closed by the constitutive law of the
material
\begin{equation}
\label{constitutive}
\bfD(t,\bfx) = \eps_0 \eps_\infty\bfE(t,\bfx) 
+ \eps_0 \int_0^t \bfE(t-\tau,\bfx) \chi(\tau)d\tau,
\end{equation}
where $\eps_\infty$ is the relative permittivity at the infinite frequency and
$\chi$ the linear susceptibility. If we discretize the integral equation
\eqref{constitutive}, we obtain what is called a recursive scheme (see
e.g. \cite{Luebbers-Hunsberger-Kunz-Standler-Schneider90},
\cite{Young-Kittichartphayak-Kwok-Sullivan95}). We can also differentiate
Eq. \eqref{constitutive} to obtain a time-differential equation for $\bfD$ 
which depends on the specific form of $\chi$. For a Debye medium, this
differential equation reads 
\begin{equation}
\label{DebyeD}
t_\rmr \partial_t \bfD + \bfD = t_\rmr \eps_0 \eps_\infty \partial_t \bfE 
+ \eps_0\eps_\rms \bfE,
\end{equation}
where $t_\rmr$ is the relaxation time and $\eps_\rms$ the static relative
permittivity. We can derive an equivalent form dealing with the polarisation 
polarisation $\bfP(t,\bfx)= \bfD(t,\bfx)-\eps_0 \eps_\infty\bfE(t,\bfx)$, 
namely
\begin{equation}
\label{DebyeP}
t_\rmr \partial_t \bfP + \bfP = \eps_0(\eps_\rms-\eps_\infty) \bfE.
\end{equation}
For a Lorentz medium with one resonant frequency $\omega_1$, we have similarly
\begin{equation}
\label{LorentzD}
\partial_t^2 \bfD + \nu \partial_t \bfD + \omega_1^2\bfD 
= \eps_0\eps_\infty\partial_t^2 \bfE + \eps_0\eps_\infty\nu \partial_t \bfE 
+ \eps_0\eps_\rms\omega_1^2\bfE,
\end{equation}
where $\nu$ is a damping coefficient and 
\begin{equation}
\label{LorentzP}
\partial_t^2 \bfP + \nu \partial_t \bfP + \omega_1^2\bfP 
= \eps_0(\eps_\rms-\eps_\infty)\omega_1^2\bfE.
\end{equation}
Denoting by $\bfJ$ the time derivative of $\bfP$, Maxwell system 
\eqref{Maxwell} can be cast as
\begin{equation}
\label{MaxwellP}
\begin{array}{rcl}
\ds\partial_t \bfB(t,\bfx) & = & \ds- \curl \bfE(t,\bfx), \\
\ds \eps_0\eps_\infty \partial_t \bfE(t,\bfx) 
& = & \ds \frac1{\mu_0} \curl \bfB(t,\bfx) - \bfJ(t,\bfx).
\end{array}
\end{equation}

\subsection{Yee Scheme}

To discretize Maxwell equations in a passive medium ($\bfJ=0$), we use Yee
scheme \cite{Yee66}, which consists in staggering space and time discretization
grids for the different fields. We denote by 
$c_\infty = 1/\sqrt{\eps_0\eps_\infty\mu_0}$ the light speed at infinite
frequency. If the space step $\dx$ is the same in all directions and $\dt$ is
the time step, the CFL condition is $c_\infty \dt/\dx\leq1$ in space dimension
$N=1$ and $c_\infty \dt/\dx\leq1/\sqrt2$ for $N=2$ or 3. In dimension 1, we can
for example only consider fields $E\equiv E_x$ et $B\equiv B_y$ which discrete
equivalents are $E_j^n \simeq E(n\dt,j\dx)$ (with similar notations for 
$D\equiv D_x$) and 
$B_{j+\frac12}^{n+\frac12} \simeq B((n+\frac12)\dt,(j+\frac12)\dx)$. Yee scheme
for the initial Maxwell system \eqref{Maxwell} in variables $\bfE$, $\bfB$ and
$\bfD$ therefore reads
\begin{equation}
\label{Max}
\begin{array}{rcl}
\ds \frac1\dt (B_{j+\frac12}^{n+\frac12}-B_{j+\frac12}^{n-\frac12})
& = & \ds - \frac1\dx (E_{j+1}^n-E_j^n), \\
\ds \frac1\dt (D_j^{n+1}-D_j^n) & = & \ds - \frac1{\mu_0 \dx}  
(B_{j+\frac12}^{n+\frac12}-B_{j-\frac12}^{n+\frac12}). 
\end{array}
\end{equation}
In the same way, for Maxwell system \eqref{MaxwellP} in variables $\bfE$, 
$\bfB$ and $\bfJ$, we have the Yee discretization
\begin{equation}
\label{MaxP}
\begin{array}{rcl}
\ds \frac1\dt (B_{j+\frac12}^{n+\frac12}-B_{j+\frac12}^{n-\frac12})
& = & \ds - \frac1\dx (E_{j+1}^n-E_j^n), \\
\ds \frac{\eps_0\eps_\infty}\dt (E_j^{n+1}-E_j^n) & = & 
\ds - \frac1{\mu_0 \dx} (B_{j+\frac12}^{n+\frac12}-B_{j-\frac12}^{n+\frac12})
- J_j^{n+\frac12}. 
\end{array}
\end{equation}

For the matter equations, we address "direct integration" schemes which
discretize the differential equations \eqref{DebyeD}--\eqref{LorentzP} (see 
\cite{Kashiwa-Yoshida-Fukai90}, \cite{Joseph-Hagness-Taflove91},
\cite{Young95}).

Before describing and analysing the schemes one by one, we give below the
principle of the von Neumann analysis which allows us to study their stability.

\section{Principles of the von Neumann Analysis}
\label{sec-vonNeumann}

\subsection{Schur and von Neumann polynomials}

We define two families of polynomials: Schur and simple von Neumann polynomials.
\begin{definition}
A polynomial is a Schur polynomial if all its roots $r$ satisfy $|r|<1$.
\end{definition}
\begin{definition}
A polynomial is a simple von Neumann polynomial if all its roots $r$ belong to
the unit disk ($|r|\leq1$) and all the roots of modulus 1 are simple roots.
\end{definition}

It may be difficult to localise roots of a polynomial with complicated
coefficients. On the other hand, we can turn this difficult problem into the
solving of many simpler small problems. To this aim, we construct a polynomial
series with strictly decreasing degree. To a polynomial $\phi$ defined by
\begin{equation*}
\phi(z) = c_0 + c_1 z + \dots + c_p z^p,
\end{equation*}
where $c_0$, $c_1$ \dots, $c_p\in\bbC$ and $c_p\neq0$, we associate its
conjugate polynomial $\phi^*$ which reads
\begin{equation*}
\phi^*(z) = c_p^* + c_{p-1}^* z + \dots + c_0^* z^p.
\end{equation*}
Given a polynomial $\phi_0$, we can define a series of polynomials by recursion
\begin{equation*}
\phi_{m+1}(z) = \frac{\phi_m^*(0)\phi_m(z)-\phi_m(0)\phi_m^*(z)}{z}.
\end{equation*}
This series is finite since it is clearly strictly degree decreasing: 
$\textrm{deg} \phi_{m+1} < \textrm{deg} \phi_m$,  if $\phi_m\not\equiv 0$. 
Besides, we have the following two theorems at our disposal.
\begin{theorem}
\label{Th_Schur}
A polynomial $\phi_m$ is a Schur polynomial of exact degree $d$ if and only if
$\phi_{m+1}$ is a Schur polynomial of exact degree $d-1$ and 
$|\phi_m(0)|<|\phi_m^*(0)|$.
\end{theorem}
\begin{theorem}
\label{Th_vonNeumann}
A polynomial $\phi_m$ is a simple von Neumann polynomial if and only if \\
\hspace*{1cm} $\phi_{m+1}$ is a simple von Neumann simple polynomial and 
$|\phi_m(0)|<|\phi_m^*(0)|$, \\
or \\
\hspace*{1cm} $\phi_{m+1}$ is identically zero and $\phi'_m$ is a Schur
polynomial.
\end{theorem}

To localise roots of $\phi_0$ in the unit disk or not, we only have to check
conditions at each step $m$ (non zero leading coefficient,
$|\phi_m(0)|<|\phi_m^*(0)|$, \dots) until we obtain a negative answer or a
polynomial of degree 1. 

The proofs of the above results are based on Rouch\'e theorem and are given in
\cite{Strikwerda89}.

\subsection{Stability Analysis}

The models we consider are linear. They can therefore be analysed in the
frequency domain. Hence we assume that the scheme deals with a variable
$U^n_\bfj$ with space dependency in the form 
\begin{equation*}
U^n_\bfj = U^n \exp(i\bfxi\cdot\bfj),
\end{equation*}
where $\bfxi$ et $\bfj\in\bbR^N$, $N=1,2,3$. Let $G$ be the matrix such that 
$U^{n+1}=GU^n$ and we assume it does not depend on time, nor on $\dx$ and $\dt$
separately but only on the ratio $\dx/\dt$. Let $\phi_0$ be the characteristic
polynomial $G$, then we have the following sufficient stability condition.
\begin{theorem}
\label{Th_CSstab}
A sufficient stability condition is that $\phi_0$ is a simple von Neumann
polynomial.
\end{theorem}
This condition is not a necessary one. The stability is linked to the fact that
$U^n=G^nU^0$ and corresponds to the boundedness of the iterates $G^n$ of the
matrix $G$. The case of multiple unit modulus roots can give rise to iterates 
of $G$ which are bounded (e.g. for the identity matrix) or not. For example
\begin{equation*}
\left(\begin{array}{cc}1 & 0 \\ 0 & 1\end{array}\right)^n
=\left(\begin{array}{cc}1 & 0 \\ 0 & 1\end{array}\right)
\textrm{ is bounded, and }
\left(\begin{array}{cc}1 & 1 \\ 0 & 1\end{array}\right)^n
=\left(\begin{array}{cc}1 & n \\ 0 & 1\end{array}\right)
\textrm{ is not bounded.}
\end{equation*}
This case occurs for the schemes we are dealing with and have to be treated
separately, without the help of the von Neumann analysis, which handles
characteristic polynomials and not the matrices they stem from, which induces a
loss of information.

\section{Debye Type Media}
\label{sec-Debye}

For Debye type media, we study two schemes. The first one is due to Joseph et
al. \cite{Joseph-Hagness-Taflove91} and consists in coupling Maxwell equations
in variables $\bfE$, $\bfB$ and $\bfD$ with the Debye model linking $\bfE$ and
$\bfD$. The second is due to Young \cite{Young95} and couples Maxwell equations
in variables $\bfE$, $\bfB$ and $\bfJ$ with the Debye model linking $\bfE$,
$\bfP$ and $\bfJ$. 


\subsection{Joseph et al. Model}

\subsubsection{Model Setting}

Maxwell system \eqref{Max} is closed by a discretization of the Debye model
\eqref{DebyeD}, namely
\begin{equation}
\label{DebyeD1}
\eps_0\eps_\infty t_\rmr \frac{E^{n+1}_j - E^n_j}\dt
+ \eps_0\eps_\rms \frac{E^{n+1}_j + E^n_j}2 
= t_\rmr \frac{D^{n+1}_j - D^n_j}\dt + \frac{D^{n+1}_j + D^n_j}2. 
\end{equation}
System \eqref{Max}--\eqref{DebyeD1} deals with the variable 
\begin{equation*}
U^n_j 
= (c_\infty B^{n-\frac12}_{j+\frac12}, E^n_j, D^n_j/\eps_0\eps_\infty)\transp 
= (\calB^{n-\frac12}_{j+\frac12}, \calE^n_j, \calD^n_j)\transp
\end{equation*}
and reads
\begin{eqnarray*}
\calB_{j+\frac12}^{n+\frac12}-\calB_{j+\frac12}^{n-\frac12}
& = & - \frac{c_\infty \dt}{\dx} (\calE_{j+1}^n-\calE_j^n), \\
\calD_j^{n+1}-\calD_j^n  & = & - \frac{c_\infty \dt}{\dx}
(\calB_{j+\frac12}^{n+\frac12}-\calB_{j-\frac12}^{n+\frac12}), \\
\calE_j^{n+1}-\calE_j^n       
+ \frac{\eps_\rms}{\eps_\infty}\frac{\dt}{2t_\rmr} (\calE_j^{n+1}+\calE_j^n) 
& = &
\calD_j^{n+1}-\calD_j^n + \frac{\dt}{2t_\rmr} (\calD_j^{n+1}+\calD_j^n).
\end{eqnarray*} 
We see that this formulation contains dimensionless parameters:
\begin{center}
\begin{tabular}{p{5cm}l}
$\lambda=c_\infty\dt/\dx$ & CFL  constant, \\ 
$\delta=\dt/2t_\rmr$ & normalised time step, \\
$\eps'_\rms=\eps_\rms/\eps_\infty$ & normalised static permittivity.
\end{tabular}
\end{center}
We write this system into the explicit form
\begin{eqnarray*}
\calB_{j+\frac12}^{n+\frac12} 
& = & \calB_{j+\frac12}^{n-\frac12} - \lambda (\calE_{j+1}^n-\calE_j^n), \\
\calD_j^{n+1} & = & \calD_j^n  
- \lambda (\calB_{j+\frac12}^{n-\frac12}-\calB_{j-\frac12}^{n-\frac12}) 
+ \lambda^2 (\calE_{j+1}^n-2\calE_j^n+\calE_{j-1}^n), \\
(1+\delta\eps'_\rms) \calE_j^{n+1} & = & (1-\delta\eps'_\rms)\calE_j^n  
+ (1+\delta)\lambda^2 (\calE_{j+1}^n-2\calE_j^n+\calE_{j-1}^n) \\ 
&& + (1+\delta)\lambda 
(\calB_{j+\frac12}^{n-\frac12}-\calB_{j-\frac12}^{n-\frac12})
+ 2\delta \calD^n_j,
\end{eqnarray*} 
which yields the amplification matrix
\begin{equation*}
G = \left(\begin{array}{ccc}
1 & -\lambda(e^{i\xi}-1) & 0 \\
- \frac{(1+\delta)\lambda(1-e^{-i\xi})}{1+\delta\eps'_\rms} 
& \frac{(1-\delta\eps'_\rms)+(1+\delta)\lambda^2(e^{i\xi}-2+e^{-i\xi})}
{1+\delta\eps'_\rms} 
& \frac{2\delta}{1+\delta\eps'_\rms} \\
-\lambda(1-e^{-i\xi}) & \lambda^2(e^{i\xi}-2+e^{-i\xi}) & 1 
\end{array}\right).
\end{equation*}
We set $\sigma=\lambda(e^{i\xi}-1)$ and 
$q = |\sigma|^2 = -\lambda^2(e^{i\xi}-2+e^{-i\xi})=4\lambda^2\sin^2(\xi/2)$.
With these notations $G$ reads
\begin{equation*}
G = \left(\begin{array}{ccc}
1 & -\sigma & 0 \\
\frac{(1+\delta)\sigma^*}{1+\delta\eps'_\rms} 
& \frac{(1-\delta\eps'_\rms)-(1+\delta)q}{1+\delta\eps'_\rms} 
& \frac{2\delta}{1+\delta\eps'_\rms} \\
\sigma^* & -q & 1 
\end{array}\right).
\end{equation*}

\subsubsection{Computation of the Characteristic Polynomial}

The characteristic polynomial of $G$ is equal to
\begin{eqnarray*}
P(Z) & = & \frac1{1+\delta\eps'_\rms}
\left|\begin{array}{ccc}
Z-1 & \sigma & 0 \\
-(1+\delta)\sigma^* 
& (1+\delta\eps'_\rms)Z-(1-\delta\eps'_\rms)+(1+\delta)q & -2\delta \\
-\sigma^* & q & Z-1 
\end{array}\right| \\
& = & \frac1{1+\delta\eps'_\rms} \Bigg( 
(Z-1)
\left|\begin{array}{cc}
(1+\delta\eps'_\rms)Z-(1-\delta\eps'_\rms)+(1+\delta) q & -2\delta \\
q & Z-1 
\end{array}\right|  
+ q \left|\begin{array}{cc}
1+\delta & -\delta \\
1 & Z-1 
\end{array}\right| \Bigg) \\
& = & \frac1{1+\delta\eps'_\rms}\Big( 
(Z-1)\left\{[1+\delta\eps'_\rms]Z^2 
- [2-(1+\delta)q] Z + [1-\delta\eps'_\rms-(1-\delta)q]\right\} \\
&& \hspace*{2cm} + q\left\{ [1+\delta] Z - [1-\delta] \right\} \Big).
\end{eqnarray*}
The characteristic polynomial is proportional to
\begin{equation*}
\phi_0(Z) = [1+\delta\eps'_\rms]Z^3 - [3 + \delta\eps'_\rms - (1+\delta)q] Z^2 
+ [3-\delta\eps'_\rms - (1-\delta)q] Z - [1-\delta\eps'_\rms].
\end{equation*}

\subsubsection{Von Neumann Analysis}

From the polynomial $\ph_0$, we perform the recursive construction of the
above-mentioned series of polynomials. We therefore define
\begin{equation*}
\ph^*_0(Z) = [1+\delta\eps'_\rms] - [3+\delta\eps'_\rms-(1+\delta)q]Z 
+ [3-\delta\eps'_\rms-(1-\delta)q]Z^2 - [1-\delta\eps'_\rms] Z^3.
\end{equation*}
The condition $|\ph_0(0)|<|\ph_0^*(0)|$ is valid. We define by recursion
\begin{eqnarray*}
\ph_1(Z) & = & \frac1Z \{\ph_0^*(0)\ph_0(Z)-\ph_0(0)\ph_0^*(Z)\} \\
& = & 2\delta \{ 2\eps'_\rms Z^2 - [4\eps'_\rms-(\eps'_\rms+1)q] Z 
+ [2\eps'_\rms-(\eps'_\rms-1)q] \}, \\
\ph_1^*(Z) & = & 2\delta \{ 2\eps'_\rms - [4\eps'_\rms-(\eps'_\rms+1)q] Z 
+ [2\eps'_\rms-(\eps'_\rms-1)q] Z^2 \}.
\end{eqnarray*}
Since $\eps_\rms \geq \eps_\infty$, we have $\eps'_\rms\geq1$ and the quantity 
$\eps'_\rms-1$ is nonnegative. If $q=0$ or $\eps'_\rms=1$, we have exactly 
$|\ph_1(0)|=|\ph_1^*(0)|$, and these specific cases have to be treated
separately (see below). In the opposite case, condition
$|\ph_1(0)|<|\ph_1^*(0)|$ reverts to $(\eps'_\rms-1)q < 4\eps'_\rms$. It is
reasonable to assume we will not obtain a better result than with the raw Yee
scheme ($\lambda \leq 1$) and therefore $q\in[0,4]$. In that case, and provided
$q\neq0$, we do have $|\ph_1(0)|<|\ph_1^*(0)|$. Moreover the degree of
polynomial $\ph_1$ is 2. Last
\begin{eqnarray*}
\ph_2(Z) & = & \frac1Z (\ph_1^*(0)\ph_1(Z)-\ph_1(0)\ph_1^*(Z)) \\
& = & 4\delta^2(\eps'_\rms-1)q 
\big[ (4\eps'_\rms-(\eps'_\rms-1)q)Z-(4\eps'_\rms-(\eps'_\rms+1)q) \big].
\end{eqnarray*}
Always in the case when $\eps'_\rms>1$ and $q\in]0,4]$, the leading coefficient
\begin{equation*}
4\eps'_\rms-(\eps'_\rms-1)q = \eps'_\rms (4-q) + q \neq0,
\end{equation*} 
and is degree of $\ph_2$ is 1. The root of $\ph_2$ is
\begin{equation*}
Z = \frac{4\eps'_\rms-(\eps'_\rms+1)q}{4\eps'_\rms-(\eps'_\rms-1)q}.
\end{equation*}
The modulus of this root is strictly lower than 1 if $q\neq4$ and therefore 
$\ph_2$ and hence $\ph_0$ are Schur polynomials thanks to 
Theorem~\ref{Th_Schur}. If $q=4$, the root of $\ph_2$ is $-1$ and $\ph_2$
and hence $\ph_0$ are simple von Neumann polynomials thanks to 
Theorem~\ref{Th_vonNeumann}. In both cases, we obtain the stability with the 
only assumption $\lambda<1$, provided we treat the above-mentioned special 
cases.

\subsubsection{Case $q=0$}

The case when $q=0$ corresponds to the characteristic polynomial 
\begin{equation*}
\ph_0(Z) = (Z-1)^2([1+\delta\eps'_\rms]Z-[1-\delta\eps'_\rms])
\end{equation*}
which is not a simple von Neumann one. We shall therefore study the 
amplification matrix directly, which is then simply
\begin{equation*}
G = \left(\begin{array}{ccc}
1 & 0 & 0 \\
0 & \ds \frac{1-\delta\eps'_\rms}{1+\delta\eps'_\rms}
& \ds \frac{2\delta}{1+\delta\eps'_\rms} \\
0 & 0 & 1
\end{array}\right).
\end{equation*}
We clearly see that the eigenvectors corresponding to the eigenvalue 1 are in
two stable eigensubspaces. The other eigenvalue has a modulus strictly lower
than 1. Iterates of this matrix are therefore bounded. This conclusion is valid
for all $\eps'_\rms\geq1$.

\subsubsection{Case $\eps'_\rms=1$}

The case when $\eps'_\rms=1$ gives rise to a breaking of condition 
$|\ph_1(0)|<|\ph_1^*(0)|$. We shall therefore study directly the nature of 
$\ph_1$ without carrying recursion over. We have
\begin{eqnarray*}
\ph_1(Z) & = & 4\delta \{ Z^2 - [2-q] Z + 1 \},
\end{eqnarray*}
which determinant is $q(q-4)$ and is therefore negative if $q\in]0,4[$. The
roots of $\phi_1$ are therefore complex conjugate, distinct and their modulus 
is 1 (their product is equal to 1). The polynomial $\phi_1$ is therefore a 
simple von Neumann one, and $\phi_0$ also.

\subsubsection{Case $q=4$}

The last case we have to treat is $\eps'_\rms=1$ and $q=4$, where $-1$ is a
double root of $\ph_1$. It is also a double root of $\ph_0$ which reads 
\begin{eqnarray*}
\ph_0(Z) & = & (Z+1)^2([1+\delta]Z - [1-\delta]).
\end{eqnarray*}
Hence we study directly the amplification matrix which reads simply
\begin{equation*}
G = \left(\begin{array}{ccc}
1 & -\sigma & 0 \\
\sigma^* & \frac{1-\delta}{1+\delta}-q & \frac{2\delta}{1+\delta} \\
\sigma^* & -q & 1 
\end{array}\right).
\end{equation*}
There is no trivial splitting in two distinct eigensubspaces. We compute the
eigenvectors associated to the eigenvalue $-1$. To this aim we solve
\begin{equation*}
(G+\Id)V = \left(\begin{array}{ccc}
2 & -\sigma & 0 \\
\sigma^* & \frac2{1+\delta}-q & \frac{2\delta}{1+\delta} \\
\sigma^* & -q & 2 
\end{array}\right)V=0
\ \Longleftrightarrow\
\left(\begin{array}{ccc}
2 & -\sigma & 0 \\
0 & 1 & -1 \\
\sigma^* & -q & 2 
\end{array}\right)V=0
\end{equation*}
and we only find one eigendirection, that of $V=(\sigma,2,2)\transp$. A minimal
two-dimensional eigensubspace is therefore associated to the eigenvalue $-1$ 
and iterates $G^n$ are linearly increasing with $n$. Hence we conclude to
instability when $q=4$ and $\eps'_\rms=1$.

\subsubsection{Synthesis for the Debye--Joseph et al. Model}

The scheme \eqref{Max}--\eqref{DebyeD1} for the one-dimensional 
Maxwell--Debye equation is stable with the condition
\begin{equation*}
\dt \leq \dx / c_\infty \textrm{ if } \eps_\rms > \eps_\infty 
\textrm{ \ \ and \ \ } 
\dt < \dx / c_\infty \textrm{ if } \eps_\rms = \eps_\infty.
\end{equation*}

We have already seen in this first example different types of arguments to
conclude to stability: the generic case ($q\in]0,4[$ and $\eps_\rms>1$) gives
rise to a Schur polynomial \textit{via} Theorem~\ref{Th_Schur}, the cases
$q\in]0,4[$ and $\eps'_\rms=1$ or $q=4$ and $\eps'_\rms>1$ to a simple von
Neumann polynomial \textit{via} Theorem~\ref{Th_vonNeumann} and last, the case
$q=0$ to a (not simple) von Neumann polynomial, but with a double eigenvalue
that operates on two stable and distinct eigensubspaces. We have also
encountered an instable case when $\eps'_s=1$ and $q=4$ which nevertheless
corresponds to a (non simple) von Neumann polynomial.


\subsection{Young Model}

\subsubsection{Model Setting}

Maxwell system \eqref{MaxP} is closed by two discretizations of Debye equation
\eqref{DebyeP}, namely
\begin{equation}
\label{DebyeP1}
t_\rmr \frac{P^{n+\frac12}_j - P^{n-\frac12}_j}\dt = 
- \frac{P^{n+\frac12}_j + P^{n-\frac12}_j}2 
+ \eps_0(\eps_\rms-\eps_\infty) E^n_j, 
\end{equation}
and
\begin{equation}
\label{DebyeP2}
t_\rmr J^{n+\frac12}_j = - P^{n+\frac12}_j 
+ \eps_0(\eps_\rms-\eps_\infty) \frac{E^{n+1}_j + E^n_j}2.
\end{equation}
Although we make use of $J^{n+\frac12}_j$ in the description of the scheme, 
this is not a genuine variable and the system 
\eqref{MaxP}--\eqref{DebyeP1}--\eqref{DebyeP2} deals with the variable
\begin{equation*}
U^n_j 
= (c_\infty B^{n-\frac12}_{j+\frac12}, E^n_j, 
P^{n-\frac12}_j/\eps_0\eps_\infty)\transp 
= (\calB^{n-\frac12}_{j+\frac12}, \calE^n_j, \calP^{n-\frac12}_j)\transp,
\end{equation*}
and reads
\begin{eqnarray*}
\calB_{j+\frac12}^{n+\frac12}-\calB_{j+\frac12}^{n-\frac12}
& = & - \lambda (\calE_{j+1}^n-\calE_j^n), \\
\calE_j^{n+1}-\calE_j^n  & = & - \lambda
(\calB_{j+\frac12}^{n+\frac12}-\calB_{j-\frac12}^{n+\frac12}) 
+ 2 \delta \calP_j^{n+\frac12} - \delta\alpha (\calE^{n+1}_j+\calE^n_j), \\
\calP_j^{n+\frac12}-\calP_j^{n-\frac12} 
& = & - \delta (\calP_j^{n+\frac12}+\calP_j^{n-\frac12}) 
+ 2\delta\alpha \calE^n_j.
\end{eqnarray*} 
In this system, apart from the notations $\lambda$, $\delta$ which we have
already defined, we have introduced the dimensionless parameter 
$\alpha=\eps'_\rms-1$ which is , as we already mentioned, a non negative
parameter. We rewrite this system in the explicit form
\begin{eqnarray*}
\calB_{j+\frac12}^{n+\frac12} 
& = & \calB_{j+\frac12}^{n-\frac12} - \lambda (\calE_{j+1}^n-\calE_j^n), \\
(1+\delta\alpha) \calE_j^{n+1} & = & (1-\delta\alpha) \calE_j^n  
- \lambda (\calB_{j+\frac12}^{n-\frac12}-\calB_{j-\frac12}^{n-\frac12}) 
+ \lambda^2 (\calE_{j+1}^n-2\calE_j^n+\calE_{j-1}^n) \\
&& + 2\delta \frac{1-\delta}{1+\delta} \calP_j^{n-\frac12} 
+ \frac{4\delta^2\alpha}{1+\delta} \calE^n_j, \\
(1+\delta) \calP_j^{n+\frac12} & = & (1-\delta)\calP_j^{n-\frac12}  
+ 2\delta\alpha \calE^n_j,
\end{eqnarray*} 
from which stems the amplification matrix
\begin{equation*}
G = \left(\begin{array}{ccc}
1 & -\sigma & 0 \\
\frac{\sigma^*}{1+\delta\alpha} 
& \frac{(1+\delta)(1-\delta\alpha)+4\delta^2\alpha - (1+\delta)q}
{(1+\delta)(1+\delta\alpha)} 
& \frac{1-\delta}{1+\delta}\frac{2\delta}{1+\delta\alpha}  \\
0 & \frac{2\delta\alpha}{1+\delta}& \frac{1-\delta}{1+\delta}
\end{array}\right).
\end{equation*}

\subsubsection{Computation of the Characteristic Polynomial}

The characteristic polynomial $G$ is 
\begin{equation*}
P(Z) = \left|\begin{array}{ccc}
Z-1 & \sigma & 0 \\
-\frac{\sigma^*}{1+\delta\alpha} 
& Z-\frac{(1+\delta)(1-\delta\alpha)+4\delta^2\alpha - (1+\delta)q}
{(1+\delta)(1+\delta\alpha)} 
& -\frac{1-\delta}{1+\delta}\frac{2\delta}{1+\delta\alpha} \\
0 & - \frac{2\delta\alpha}{1+\delta} & Z-\frac{1-\delta}{1+\delta} 
\end{array}\right|. 
\end{equation*}
To reduce computations, we set $Y=Z-1$, which yields \\
$(1+\delta)^2(1+\delta\alpha) P(Z)$
\begin{eqnarray*}
& = &
\left|\begin{array}{ccc}
Y & \sigma & 0 \\
-(1+\delta)\sigma^*
& (1+\delta)(1+\delta\alpha)Y+2\delta\alpha(1-\delta)+(1+\delta)q 
& -2\delta(1-\delta) \\
0 & -2\delta\alpha & (1+\delta)Y+2\delta 
\end{array}\right| \\
& = &
\left|\begin{array}{ccc}
Y & \sigma & 0 \\
-(1+\delta)\sigma^* & (1+\delta)(1+\delta\alpha)Y+(1+\delta)q 
& (1+\delta)(1-\delta)Y \\
0 & -2\delta\alpha & (1+\delta)Y+2\delta 
\end{array}\right|.
\end{eqnarray*}
We see that $(1+\delta)$ is a factor in both sides and therefore
\begin{eqnarray*}
(1+\delta)(1+\delta\alpha) P(Z) & = & Y
\left|\begin{array}{cc}
(1+\delta\alpha)Y+q & (1-\delta)Y \\
-2\delta\alpha & (1+\delta)Y+2\delta 
\end{array}\right| 
+ q
\left|\begin{array}{cc}
1 & -2\delta(1-\delta) \\
0 & (1+\delta)Y+2\delta 
\end{array}\right| \Bigg) \\
& = & 
Y\left\{[(1+\delta)(1+\delta\alpha)]Y^2 
+ [2\delta(1+\alpha)+(1+\delta)q]Y 
+ [2\delta q]\right\}\\
&& + q(1+\delta)\left\{ [1+\delta] Y + [2\delta] \right\} \\
& = & [(1+\delta)(1+\delta\alpha)]Y^3 + [2\delta(1+\alpha)+(1+\delta)q] Y^2 \\
&&+ [(1+3\delta)q] Y + [2\delta q]. 
\end{eqnarray*}
The characteristic polynomial is proportional to
\begin{eqnarray*}
\phi_0(Z) & = & 
[(1+\delta\alpha)(1+\delta)] Z^3 
- [3 + \delta + \delta\alpha + 3\delta^2\alpha - (1+\delta)q] Z^2 \\
&& + [3 - \delta - \delta\alpha + 3\delta^2\alpha - (1-\delta)q] Z 
- [(1-\delta\alpha)(1-\delta)].
\end{eqnarray*}

\subsubsection{Von Neumann Analysis}

Condition $|\ph_0(0)|<|\ph_0^*(0)|$ is valid without any assumption. We define
by recursion
\begin{eqnarray*}
\ph_1(Z) & = & 2\delta\{[2(1+\alpha)(1+\delta^2\alpha)]Z^2
-[4(1+\alpha)(1+\delta^2\alpha) - (2+\alpha+\delta^2\alpha)q] Z \\
&& + [2(1+\alpha)(1+\delta^2\alpha) - \alpha(1-\delta^2)q]\}.
\end{eqnarray*}
The case when $\delta^2>1$ does not allow to fulfil the condition 
$|\ph_1(0)|<|\ph_1^*(0)|$. We will assume therefore for the von Neumann 
analysis that $\delta<1$, which bounds the time step with respect to the time 
delay $t_\rmr$. This is reasonable from the point of view of modelling: we 
cannot approximate the delay equation with too large a time step. Such an 
assumption was however not necessary for the Joseph et al. scheme.\\
The equality case $|\ph_1(0)|=|\ph_1^*(0)|$ is obtained when $q=0$, $\alpha=0$
(i.e. $\eps'_\rms=1$) or $\delta=1$. These cases shall be treated separately
again.\\
If $\alpha>0$, $q>0$ and $\delta<1$, then $|\ph_1(0)|<|\ph_1^*(0)|$ is
equivalent to $\alpha(1-\delta^2)q < 4 (1+\alpha)(1+\delta^2\alpha)$, which is
clearly true if $q\in]0,4]$. Besides, the degree of polynomial $\ph_1$ is 2. \\
In the general case ($\alpha>0$, $q\in]0,4]$ and $\delta<1$), we then compute
$\ph_2$
\begin{eqnarray*}
\ph_2(Z) & = & 4\delta^2\alpha(1-\delta^2)q \{
[4(1+\alpha)(1+\delta^2\alpha)-\alpha(1-\delta^2)q]Z \\
&& \hspace{3cm}- [4(1+\alpha)(1+\delta^2\alpha) - (2+\alpha+\delta^2\alpha)q]
\}
\big].
\end{eqnarray*}
We split the study according to the sign of $\ph_2(0)$ ($\ph_2^*(0)$ is clearly
always positive) and in both cases $|\ph_2(0)|<|\ph_2^*(0)|$, for
$q\in]0,4]$. The root of $\ph_2$ therefore belongs to the interval $]-1,1[$ and
$\ph_0$ is a Schur polynomial. Hence we obtain the stability with the
assumptions $\lambda<1$ and $\delta<1$, provided we treat the above-mentioned
specific cases.

\subsubsection{Case $q=0$}

The case when $q=0$ corresponds to the characteristic polynomial
\begin{equation*}
\ph_0(Z) = (Z-1)^2
(Z-\frac{(1-\delta)(1-\delta\alpha)}{(1+\delta)(1+\delta\alpha)})
\end{equation*}
and is not a simple von Neumann one. The amplification matrix reads
\begin{equation*}
G = \left(\begin{array}{ccc}
1 & 0 & 0 \\
0 & \frac{(1+\delta)(1-\delta\alpha)+4\delta^2\alpha}
{(1+\delta)(1+\delta\alpha)} 
& \frac{1-\delta}{1+\delta}\frac{2\delta}{1+\delta\alpha}  \\
0 & \frac{2\delta\alpha}{1+\delta}& \frac{1-\delta}{1+\delta}
\end{array}\right).
\end{equation*}
We clearly see that the eigenvectors corresponding to the eigenvalue 1 are in
two stable eigensubspaces. Iterates of this matrix are therefore bounded. This
conclusion is once more valid in the limit cases $\delta=1$ and $\alpha=0$.

\subsubsection{Case $\eps'_\rms=1$}

We notice that $\ph_0$ is the same as that for the Joseph et al. model for
$\eps'_\rms=1$. Polynomial $\ph_0$ is therefore a simple von Neumann polynomial
for $q\neq4$ (see above). The value of $\delta$ does not play any r\^ole
here. This corresponds to different amplification matrices, operating on
different sets of variables, the link between both formulations being not
straightforward. Case $q=4$ has therefore to be treated anew.

\subsubsection{Case $q=4$}

If $q=4$, only the case when $\alpha=0$ has not been treated by the general
study and  $-1$ is once more a double eigenvalue
\begin{equation*}
(G+\Id)V = \left(\begin{array}{ccc}
2 & -\sigma & 0 \\
\sigma^* & 2-q & 2\delta \frac{1-\delta}{1+\delta} \\
0 & 0 & \frac2{1+\delta}
\end{array}\right)V=0 
\ \Longleftrightarrow\ 
\left(\begin{array}{ccc}
2 & -\sigma & 0 \\
\sigma^* & 2-q & 0 \\
0 & 0 & \frac2{1+\delta}
\end{array}\right)V=0, 
\end{equation*}
and the only eigendirection is that of $V=(\sigma,2,0)\transp$, which gives 
rise to linearly increasing iterates $G^n$ and to instabilities. The fact that
$\delta=1$ or not does not play any r\^ole in this argument.

\subsubsection{Case $\delta=1$}

There remains to study the case $\delta=1$ for $q\in]0,4]$ and $\alpha>0$.
Then $Z=0$ is a trivial root of $\ph_0$, which simply reads
\begin{eqnarray*}
\ph_0(Z)= 2(1+\alpha)Z\{Z^2-[2-\frac{q}{1+\alpha}]Z+1\}.
\end{eqnarray*}
The discriminant of the second order factor is
$\Delta=\frac{q}{(1+\alpha)^2}(q-4(1+\alpha))<0$. We therefore have two 
distinct complex conjugate eigenvalues of modulus 1. Polynomial $\ph_0$ is a 
simple von Neumann polynomial.

\subsubsection{Synthesis for the Debye--Young Model}

The scheme \eqref{MaxP}--\eqref{DebyeP1}--\eqref{DebyeP2} for the
one-dimensional Maxwell--Debye equation is stable with the condition
\begin{equation*}
\dt \leq \min (\dx / c_\infty ,2 t_\rmr ) \textrm{ if } \eps_\rms > \eps_\infty
\textrm{ \ \ and \ \ } 
\dt < \dx / c_\infty \textrm{ if } \eps_\rms = \eps_\infty.
\end{equation*}

If $\eps_\rms>\eps_\infty$, the stability condition is more restrictive for the
Young scheme than for Joseph et al. scheme. The obtained bound is also related
to the good approximation of Debye equation.

\section{Anharmonic Lorentz Type Media}
\label{sec-Lorentz}

For Lorentz type media, we study three schemes. The first one is due to Joseph
et al. \cite{Joseph-Hagness-Taflove91} and consists in coupling Maxwell
equations in the variables $\bfE$, $\bfB$ and $\bfD$ with the Lorentz model
linking $\bfE$ and $\bfD$. The second and third are due to Kashiwa et
al. \cite{Kashiwa-Yoshida-Fukai90} and Young \cite{Young95} respectively and
both couple Maxwell equations in the variables $\bfE$, $\bfB$ and $\bfJ$ with
the Lorentz model linking $\bfE$, $\bfP$ and $\bfJ$. They differ in the choice
of the time discretization of $\bfJ$. \\

We restrict here to the anharmonic case for which the damping $\nu$ is
non-zero. The harmonic case ($\nu=0$) is treated with the same schemes but the
analysis happens to be much more technical. To keep proofs readable in the 
general case we postpone the harmonic case to the next section.


\subsection{Joseph et al. Model}

\subsubsection{Model Setting}

Maxwell system \eqref{Max} is closed by a discretization of the Lorentz 
equation \eqref{LorentzD}, namely 
\begin{equation}
\label{LorentzD1}
\begin{array}{l}
\ds \eps_0\eps_\infty \frac{E^{n+1}_j - 2E^n_j + E^{n-1}_j}{\dt^2}
+ \nu \eps_0\eps_\infty \frac{E^{n+1}_j - E^{n-1}_j}{2\dt} 
+ \eps_0\eps_\rms \omega_1^2 \frac{E^{n+1}_j + E^{n-1}_j}2 \\ 
\ds \hspace{2cm} 
= \frac{D^{n+1}_j - 2D^n_j + D^{n-1}_j}{\dt^2}
+ \nu \frac{D^{n+1}_j - D^{n-1}_j}{2\dt} 
+ \omega_1^2 \frac{D^{n+1}_j + D^{n-1}_j}2.
\end{array}
\end{equation}
The explicit version of system \eqref{Max}--\eqref{LorentzD1} does not use
explicitly the variable $D^{n-1}_j$. Indeed we can use the explicit formula to
compute $D_j^{n+1}-D_j^n$ and the implicit one to compute $D_j^n-D_j^{n-1}$ and
therefore the system deals with the variable
\begin{equation*}
U^n_j 
= (c_\infty B^{n-\frac12}_{j+\frac12}, E^n_j, E^{n-1}_j, D^n_j/\eps_0\eps_\infty)\transp 
= (\calB^{n-\frac12}_{j+\frac12}, \calE^n_j, \calE^{n-1}_j, \calD^n_j)\transp
\end{equation*}
and reads
\begin{eqnarray*}
\calB_{j+\frac12}^{n+\frac12} 
& = & \calB_{j+\frac12}^{n-\frac12} - \lambda (\calE_{j+1}^n-\calE_j^n), \\
\calD_j^{n+1} & = & \calD_j^n  
- \lambda (\calB_{j+\frac12}^{n-\frac12}-\calB_{j-\frac12}^{n-\frac12}) 
+ \lambda^2 (\calE_{j+1}^n-2\calE_j^n+\calE_{j-1}^n), \\
(1+\delta + \omega\eps'_\rms) \calE_j^{n+1}
& = & 2\calE_j^n + (1+\delta + \omega) \lambda^2 
(\calE_{j+1}^n-2\calE_j^n+\calE_{j-1}^n) 
- (1-\delta+\omega\eps'_\rms)\calE_j^{n-1} \\
& & - 2\delta\lambda (\calB_{j+\frac12}^{n-\frac12} 
- \calB_{j-\frac12}^{n-\frac12}) + 2\omega \calD_j^n.
\end{eqnarray*} 
In this system, apart from the already used notations $\lambda$ and 
$\eps'_\rms$, we have denoted
\begin{center}
\begin{tabular}{p{5cm}l}
$\delta=\dt\nu/2$ & normalised time step, \\
$\omega=\omega_1^2\dt^2/2$ & square of the normalised frequency.
\end{tabular}
\end{center}
The amplification matrix of the system is
\begin{equation*}
G = \left(\begin{array}{cccc}
1 & -\sigma & 0 & 0 \\
\frac{2\delta\sigma^*}{1+\delta+\omega\eps'_\rms} 
& \frac{2-q(1+\delta+\omega)}{1+\delta+\omega\eps'_\rms} 
& -\frac{1-\delta+\omega\eps'_\rms}{1+\delta+\omega\eps'_\rms} 
& \frac{2\omega}{1+\delta+\omega\eps'_\rms} \\
0 & 1 & 0 & 0 \\
\sigma^* & -q & 0 & 1
\end{array}\right).
\end{equation*}

\subsubsection{Computation of the Characteristic Polynomial}

The characteristic polynomial of $G$ is equal to 
\begin{equation*}
P(Z) = \left|\begin{array}{cccc}
Z-1 & \sigma & 0 & 0 \\
-\frac{2\delta\sigma^*}{1+\delta+\omega\eps'_\rms} 
& Z-\frac{2-q(1+\delta+\omega)}{1+\delta+\omega\eps'_\rms} 
& \frac{1-\delta+\omega\eps'_\rms}{1+\delta+\omega\eps'_\rms} 
& -\frac{2\omega}{1+\delta+\omega\eps'_\rms} \\
0 & -1 & Z & 0 \\
-\sigma^* & q & 0 & Z-1
\end{array}\right|. 
\end{equation*}
Therefore \\
$(1+\delta+\omega\eps'_\rms) P(Z)$
\begin{eqnarray*}
& = & \left|\begin{array}{cccc}
Z-1 & \sigma & 0 & 0 \\
-2\delta\sigma^* 
& (1+\delta+\omega\eps'_\rms)Z-2+(1+\delta+\omega)q
& 1-\delta+\omega\eps'_\rms & -2\omega \\
0 & -1 & Z & 0 \\
-\sigma^* & q & 0 & Z-1
\end{array}\right| \\
& = & \left|\begin{array}{ccc}
Z-1 & 0 & 0 \\
-2\delta\sigma^* 
& 1-\delta+\omega\eps'_\rms & -2\omega \\
-\sigma^* & 0 & Z-1
\end{array}\right|  
+ Z \left|\begin{array}{ccc}
Z-1 & \sigma & 0 \\
- 2\delta\sigma^* 
& (1+\delta+\omega\eps'_\rms)Z-2+(1+\delta+\omega)q & -2\omega \\
-\sigma^* & q & Z-1
\end{array}\right| \\
& = & (Z-1)^2[1-\delta+\omega\eps'_\rms] \\
&& + Z(Z-1) 
\{ [(1+\delta+\omega\eps'_\rms)Z-2+(1+\delta+\omega)q][Z-1]+2\omega q \}  
+ qZ[2\delta(Z-1)+2\omega].
\end{eqnarray*}
The characteristic polynomial is proportional to
\begin{eqnarray*}
\phi_0(Z) & = & 
[1+\delta+\omega\eps'_\rms] Z^4 
- [4 + 2\delta + 2\omega\eps'_\rms - (1 + \delta + \omega) q] Z^3 
+ [6 + 2\omega\eps'_\rms - 2q] Z^2 \\
&& - [4 - 2\delta + 2\omega\eps'_\rms - (1-\delta+\omega) q] Z 
+ [1-\delta+\omega\eps'_\rms].
\end{eqnarray*}

\subsubsection{Von Neumann Analysis}

We successively compute
\begin{eqnarray*}
\ph_1(Z) & = & 2\delta \{ [2(1+\omega\eps'_\rms)] Z^3 
- [6+4\omega\eps'_\rms - (2+\omega(1+\eps'_\rms))q] Z^2 \\
&& + [6+2\omega\eps'_\rms -2q] Z - [2+\omega(\eps'_\rms-1)q] \}, 
\allowdisplaybreaks \\
\ph_2(Z) & = & 4\delta^2\omega \{ [4\eps'_\rms(2+\omega\eps'_\rms) 
- 4(\eps'_\rms-1)q - \omega (\eps'_\rms-1)^2 q^2] Z^2 \\
&& -[8\eps'_\rms(2+\omega\eps'_\rms) 
- 4((\eps'_\rms-1) -\eps'_\rms(2+\omega\eps'_\rms))q 
+ 2(\eps'_\rms-1)q^2] Z \\
&& + [4\eps'_\rms(2+\omega\eps'_\rms) - 4 (\eps'_\rms-1)(2+\omega\eps'_\rms) q
+(2+ \omega(1+\eps'_\rms)) (\eps'_\rms-1) q^2]
\}, \allowdisplaybreaks \\
\ph_3(Z) & = & 64 \delta^4\omega^2 (\eps'_\rms-1)(1+\omega\eps'_\rms)q (2-q)
\times \\
& \times & \{
[4\eps'_\rms(2+\omega\eps'_\rms)-(\eps'_\rms-1)(6+2\omega\eps'_\rms))q 
+ (\eps'_\rms-1)(1+\omega)q^2] Z \\
&& \hspace*{1cm}
-[4\eps'_\rms(2+\omega\eps'_\rms)-2((\eps'_\rms-1)
+\eps'_\rms(2+\omega\eps'_\rms))q + (\eps'_\rms-1)q^2]
\}.
\end{eqnarray*}
The root of $\ph_3$ is
\begin{equation*}
Z = \frac{(2-q)(2\eps'_\rms(2+\omega\eps'_\rms)-(\eps'_\rms-1)q)}
{(2-q)(2\eps'_\rms(2+\omega\eps'_\rms)-(\eps'-1)q)+2(2+\omega\eps'_\rms)q 
+ (\eps'_\rms-1)(1+\omega)q^2},
\end{equation*}
form on which we easily see that $Z$ remains of modulus $<1$ if $q\in]0,2[$,
which corresponds to the condition for a multi-dimensional Yee scheme 
($\lambda\leq1/\sqrt2$). We notice from now on that we shall treat cases $q=0$
and $q=2$ apart because $\ph_3\equiv0$. In the general case, we have to check
the intermediate properties. First, the degree of polynomials $\ph_1$ and
$\ph_2$ is 3 and 2 respectively. The degree of polynomial $\ph_3$ is 1 provided
$\eps'_\rms>1$. We shall treat the case $\eps'_\rms=1$ apart. In the general
case ($q\neq0$ and $\eps'_\rms>1$), there remains to check the estimates 
between
\begin{eqnarray*}
\ph_0(0) & = & 1-\delta+\omega\eps'_\rms, \\
\ph^*_0(0) & = & 1+\delta+\omega\eps'_\rms, \\
\ph_1(0) & = & 2\delta\big[ -2-\omega(\eps'_\rms-1)q \big], \\
\ph^*_1(0) & = & 2\delta\big[ 2(1+\omega\eps'_\rms) \big], 
\allowdisplaybreaks \\
\ph_2(0) & = & 4\delta^2\omega \big[
4\eps'_\rms(2+\omega\eps'_\rms) - 4 (\eps'_\rms-1)(2+\omega\eps'_\rms) q 
+(2+ \omega(1+\eps'_\rms)) (\eps'_\rms-1) q^2
\big], \\
\ph^*_2(0) & = &  4\delta^2\omega \big[ 
4\eps'_\rms(2+\omega\eps'_\rms) - 4(\eps'_\rms-1)q 
- \omega (\eps'_\rms-1)^2 q^2 \big].
\end{eqnarray*}
It is clear that for $q\in]0,2[$, we have $|\ph_0(0)|<|\ph^*_0(0)|$ and 
$|\ph_1(0)|<|\ph^*_1(0)|$. A simple calculation shows that 
\begin{eqnarray*}
\ph_2(0) & = & 4 \delta^2\omega \big[
(2-q)^2 (2(1+\omega)(\eps'_\rms-1)+\omega(\eps'_\rms-1)^2) 
+ 8 + 4  \omega + 4\omega(\eps'_\rms-1) q
\big], \\
\ph^*_2(0) & = & 4 \delta^2\omega \big[
(2-q)(4(\eps'_\rms-1)+\omega(\eps'_\rms-1)^2q) 
+ 8 + 4  \omega + 8\omega(\eps'_\rms-1)
\big],
\end{eqnarray*}
form on which we readily see that both quantities are positive. Besides
\begin{equation*}
\ph^*_2(0)-\ph_2(0) = 8 \delta^2\omega (1+\omega\eps'_\rms)(\eps'_\rms-1)
q (2-q) > 0. 
\end{equation*}
We therefore checked all the assumptions. In the general case, $\ph_0$ is a 
Schur polynomial.

\subsubsection{Case $q=0$}

The case $q=0$ gives anew rise to a separate study. We have
\begin{eqnarray*}
\ph_0(Z) & = & (Z-1)^2 \big[ (Z-1)^2 + \delta (Z^2-1) 
+ \omega\eps'_\rms (Z^2+1)\big].
\end{eqnarray*} 
The corresponding amplification matrix is
\begin{equation*}
\left(\begin{array}{cccc}
1 & 0 & 0 & 0 \\
0 & \ds \frac{2}{1+\delta+\omega\eps'_\rms} 
& \ds - \frac{1-\delta+\omega\eps'_\rms}{1+\delta+\omega\eps'_\rms} 
& \ds \frac{2\omega}{1+\delta+\omega\eps'_\rms} \\
0 & 1 & 0 & 0 \\
0 & 0 & 0 & 1
\end{array}\right).
\end{equation*}
Once more, 1 is a double root but in two distinct eigensubspaces. There 
remains to check that the other factor of the polynomial, namely 
\begin{equation*}
\psi_0(Z) = [1+\delta+\omega\eps'_\rms]Z^2 - 2 Z + [1-\delta+\omega\eps'_\rms],
\end{equation*}
is a Schur (or a simple von Neumann) one. We do have
$|\psi_0(0)| < |\psi^*_0(0)|$ and we compute
\begin{equation*}
\psi_1(Z) = 4\delta\{[1+\omega\eps'_\rms]Z-1\}.
\end{equation*}
The modulus of both remaining eigenvalues is strictly less that 1 and iterates
of the amplification matrix are bounded. This holds whatever the value of
$\eps'_\rms$.

\subsubsection{Case $q=2$}

In the specific case $q=2$, $\ph_0$ reads
\begin{eqnarray*}
\phi_0(Z) & = &
[1+\delta+\omega\eps'_\rms] Z^4 
- [2 + 2\omega(\eps'_\rms-1)] Z^3 
+ [2 + 2\omega\eps'_\rms] Z^2 
- [2 + 2\omega(\eps'_\rms-1)] Z 
+ [1-\delta+\omega\eps'_\rms] \\
& = & (Z^2+1)\{[1+\delta+\omega\eps'_\rms] Z^2 
- [2 + 2\omega(\eps'_\rms-1)] Z + [1-\delta+\omega\eps'_\rms]\},
\end{eqnarray*}
which has $\pm i$ as simple roots. This is therefore a good candidate to be a
simple von Neumann polynomial. The remains to study the other factor of the 
polynomial
\begin{equation*}
\psi_0(Z) = [1+\delta+\omega\eps'_\rms] Z^2 
- [2 + 2\omega(\eps'_\rms-1)] Z + [1-\delta+\omega\eps'_\rms]
\end{equation*}
which has not $\pm i$ as roots. We notice that $|\psi_0(0)|<|\psi_0^*(0)|$ and 
compute 
\begin{equation*}
\psi_1(Z) = 4\delta\{[1+\omega\eps'_\rms] Z - [1 + \omega(\eps'_\rms-1)]\},
\end{equation*}
which is a Schur polynomial for all $\eps'_\rms$. Polynomial $\psi_0$ is 
therefore a Schur polynomial and $\ph_0$ is a simple von Neumann polynomial. 

\subsubsection{Case $\eps'_\rms=1$}

If $\eps'_\rms=1$, the polynomial $\ph_3$ is identically zero and 
\begin{eqnarray*}
\ph_2(Z) & = & 16 \delta^2 \omega (2+\omega) 
\{Z^2 - [2-q] Z + 1 \}, \\
\ph'_2(Z) & = & 16 \delta^2 \omega (2+\omega) \{ 2Z - [2-q] \}.
\end{eqnarray*}
The root of $\ph'_2$ does have a $<1$ modulus if $q\in]0,2]$. The polynomial 
$\ph_0$ is therefore a simple von Neumann polynomial. 

\subsubsection{Synthesis for the Lorentz--Joseph et al. Model}

The scheme \eqref{Max}--\eqref{LorentzD1} for the one-dimensional anharmonic
Maxwell--Lorentz equations is stable with the condition
\begin{equation*}
\dt \leq \dx / \sqrt2 c_\infty.
\end{equation*}


\subsection{Kashiwa et al. Model}

\subsubsection{Model Setting}

A modified version of Maxwell system \eqref{MaxP} is closed by a discretization
of Lorentz equation \eqref{LorentzP}, namely
\begin{equation}
\label{LorentzP1}
\begin{array}{rcl}
\ds \frac1\dt (B_{j+\frac12}^{n+\frac12}-B_{j+\frac12}^{n-\frac12})
& = & \ds - \frac1\dx (E_{j+1}^n-E_j^n), \\
\ds \frac{\eps_0\eps_\infty}\dt (E_j^{n+1}-E_j^n)  
& = & \ds - \frac1{\mu_0 \dx} 
(B_{j+\frac12}^{n+\frac12}-B_{j-\frac12}^{n+\frac12})
- \frac1\dt (P^{n+1}_j-P^n_j), \\
\ds \frac1\dt (P^{n+1}_j-P^n_j) & = & \ds \frac12 (J^{n+1}_j+J^n_j), \\
\ds \frac1\dt (J^{n+1}_j-J^n_j) & = & \ds - \frac\nu2 (J^{n+1}_j+J^n_j) 
+ \frac{\omega_1^2(\eps_\rms-\eps_\infty)\eps_0}2 (E^{n+1}_j+E^n_j) 
- \frac{\omega_1^2}2 (P^{n+1}_j+P^n_j).
\end{array}
\end{equation}
The system \eqref{LorentzP1} deals with the variable
\begin{equation*}
U^n_j = (c_\infty B^{n-\frac12}_{j+\frac12}, E^n_j, 
P^n_j/\eps_0\eps_\infty,\dt J^n_j/\eps_0\eps_\infty)\transp
= (\calB^{n-\frac12}_{j+\frac12}, \calE^n_j,\calP^n_j,\calJ^n_j)\transp
\end{equation*}
and reads
\begin{eqnarray*}
\calB_{j+\frac12}^{n+\frac12} 
& = & \calB_{j+\frac12}^{n-\frac12} - \lambda (\calE_{j+1}^n-\calE_j^n), 
\allowdisplaybreaks \\
{}[1+\delta+\frac12\omega\eps'_\rms] \calE^{n+1}_j & = & 
[1+\delta-\frac12\omega(\eps'_\rms-2)]\calE^n_j 
+ \lambda^2(1+\delta+\frac12\omega)(\calE^n_{j+1}-2\calE^n_j+\calE^n_{j-1}) \\ 
&&-\lambda (1+\delta+\frac12\omega) (\calB_{j+\frac12}^{n-\frac12} 
- \calB_{j-\frac12}^{n-\frac12})
+ \omega \calP^n_j - \calJ^n_j, 
\allowdisplaybreaks \\
{}[1+\delta+\frac12\omega\eps'_\rms] \calP^{n+1}_j & = & 
[1+\delta+\frac12\omega(\eps'_\rms-2)]\calP^n_j 
-\frac12\lambda\omega(\eps'_\rms-1) (\calB_{j+\frac12}^{n-\frac12} 
- \calB_{j-\frac12}^{n-\frac12}) \\
&& + \omega(\eps'_\rms-1) \calE^n_j 
+ \frac12 \lambda^2\omega(\eps'_\rms-1)(\calE^n_{j+1}-2\calE^n_j+\calE^n_{j-1})
+ \calJ^n_j, 
\allowdisplaybreaks \\
{}[1+\delta+\frac12\omega\eps'_\rms] \calJ^{n+1}_j & = & 
[1-\delta-\frac12\omega\eps'_\rms]\calJ^n_j 
-\lambda\omega(\eps'_\rms-1) (\calB_{j+\frac12}^{n-\frac12} 
- \calB_{j-\frac12}^{n-\frac12}) \\
&& + 2\omega(\eps'_\rms-1) \calE^n_j 
+ \lambda^2\omega(\eps'_\rms-1)(\calE^n_{j+1}-2\calE^n_j+\calE^n_{j-1}) 
- 2\omega \calP^n_j,
\end{eqnarray*} 
from which stems the amplification matrix
\begin{equation*}
G = \left(\begin{array}{cccc}
1 & -\sigma & 0 & 0 \\
\frac{\sigma^*(D-\frac12\omega(\eps'_\rms-1))}{D} 
& \frac{(1-q)D-(2-q)\frac12\omega(\eps'_\rms-1)}{D} 
& \frac{\omega}{D} & \frac{-1}{D} \\
\frac{\sigma^*\frac12\omega(\eps'_\rms-1)}{D} 
& \frac{(2-q)\frac12\omega(\eps'_\rms-1)}{D} 
& \frac{D-\omega}{D} & \frac{1}{D} \\
\frac{\sigma^*\omega(\eps'_\rms-1)}{D} 
& \frac{(2-q)\omega(\eps'_\rms-1)}{D} & \frac{-2\omega}{D} 
& \frac{2-D}{D}
\end{array}\right)
\end{equation*}
where, along with earlier notations, $D=1+\delta+\frac12\omega\eps'_\rms$.

\subsubsection{Computation of the Characteristic Polynomial}

The characteristic polynomial of $G$ is 
\begin{equation*}
P(Z) = \left|\begin{array}{cccc}
Z-1 & \sigma & 0 & 0 \\
\frac{-\sigma^*(D-\frac12\omega(\eps'_\rms-1))}{D} 
& Z-\frac{(1-q)D-(2-q)\frac12\omega(\eps'_\rms-1)}{D} 
& -\frac{\omega}{D} & \frac{1}{D} \\
\frac{-\sigma^*\frac12\omega(\eps'_\rms-1)}{D} 
& -\frac{(2-q)\frac12\omega(\eps'_\rms-1)}{D} 
& Z-\frac{D-\omega}{D} & \frac{-1}{D} \\
\frac{-\sigma^*\omega(\eps'_\rms-1)}{D} 
& -\frac{(2-q)\omega(\eps'_\rms-1)}{D} & \frac{2\omega}{D} 
& Z-\frac{2-D}{D}
\end{array}\right|. 
\end{equation*}
hence setting $X=D(Z-1)$ \\
\begin{eqnarray*}
D^4P(Z) & = & \left|\begin{array}{cccc}
X & D\sigma & 0 & 0 \\
-\sigma^*(D-\frac12\omega(\eps'_\rms-1))
& X+ qD+(2-q)\frac12\omega(\eps'_\rms-1) & -\omega & 1 \\
-\sigma^*\frac12\omega(\eps'_\rms-1) & -(2-q)\frac12\omega(\eps'_\rms-1)
& X+\omega & -1 \\
-\sigma^*\omega(\eps'_\rms-1) & -(2-q)\omega(\eps'_\rms-1) & 2\omega & X+2(D-1)
\end{array}\right| \\
& = & \left|\begin{array}{cccc}
X & D\sigma & 0 & 0 \\
-\sigma^*D
& X + qD & X & 0 \\
-\sigma^*\frac12\omega(\eps'_\rms-1) & -(2-q)\frac12\omega(\eps'_\rms-1)
& X+\omega & -1 \\
0 & 0 & -2X & X+2D
\end{array}\right| \allowdisplaybreaks \\
& = & X \left|\begin{array}{ccc}
X + qD & X & 0 \\
-(2-q)\frac12\omega(\eps'_\rms-1) & X+\omega & -1 \\
0 & -2X & X+2D
\end{array}\right|
+ qD \left|\begin{array}{ccc}
D & X & 0 \\
\frac12\omega(\eps'_\rms-1) & X+\omega & -1 \\
0 & -2X & X+2D
\end{array}\right| \allowdisplaybreaks \\
& = & X\{(X+qD)(X+\omega)(X+2D)
+(2-q)\frac12\omega(\eps'_\rms-1)X(X+2D)-2X(X+qD)\}\\
&& + qD\{D(X+\omega)(X+2D)-2DX -\frac12\omega(\eps'_\rms-1) X(X+2D)\} 
\allowdisplaybreaks \\
& = & X^4 + [2D - 2 + \omega\eps'_\rms + (D-\frac12\omega(\eps'_\rms-1))q]X^3
+ D[2\omega\eps'_\rms + (3D -2 + \frac52\omega-\frac32\omega\eps'_\rms)q]X^2
\\
&& + D^2[(2D-2+4\omega-\omega\eps'_\rms)q] X + D^3[2\omega q] 
\allowdisplaybreaks \\
& = & X^4 + [2(\delta+\omega\eps'_\rms) + (1+\delta+\frac12\omega)q]X^3 
+ D[2\omega\eps'_\rms + (1+3\delta+\frac52\omega)q]X^2 \\
&& + D^2[(2\delta+4\omega)q] X + D^3[2\omega q].
\end{eqnarray*}
The characteristic polynomial is proportional to
\begin{eqnarray*}
\phi_0(Z) & = & 
[1+\delta+\frac12\omega\eps'_\rms] Z^4
- [4+2\delta-(1+\delta+\frac12\omega)q] Z^3
+ [6-\omega\eps'_\rms+(\omega-2)q] Z^2 \\
&& - [4-2\delta-(1-\delta+\frac12\omega)q] Z 
+ [1-\delta+\frac12\omega\eps'_\rms].
\end{eqnarray*}

\subsubsection{Von Neumann Analysis}

We successively compute
\begin{eqnarray*}
\ph_1(Z) & = & 2\delta \{ [2+\omega\eps'_\rms] Z^3 
- [6+\omega\eps'_\rms - (2+\frac12\omega(\eps'_\rms+1))q] Z^2 
+ [6-\omega\eps'_\rms -(2-\omega)q] Z \\
&& \hspace{1cm}- [2-\omega\eps'_\rms+\frac12\omega(\eps'_\rms-1)q] \}, \\
\ph_2(Z) & = & 4\delta^2\omega \{
[8\eps'_\rms - (\eps'_\rms-1)(2-\omega\eps'_\rms)q
-\frac14\omega(\eps'_\rms-1)^2q^2] Z^2 \\
&&  \hspace{1cm}- [16\eps'_\rms - 8\eps'_\rms q
  + (\eps'_\rms-1)(1-\frac12\omega)q^2] Z \\
&& \hspace{1cm} + [8\eps'_\rms - (\eps'_\rms-1)(6+\omega\eps'_\rms)q
-(\eps'_\rms-1)^2(1-\frac14\omega)q^2]  \}, \\
\ph_3(Z) 
& = & 4\delta^4\omega^2 (\eps'_\rms-1)q(4-q)(2+\omega\eps'_\rms) \times \\
&& \times \{
[32\eps'_\rms - 16(\eps'_\rms-1)q + (\eps'_\rms-1)(2+\omega)q^2]Z 
-[32\eps'_\rms - 16\eps'_\rms q + (\eps'_\rms-1)(2-\omega)q^2]
\}.
\end{eqnarray*}
We see that the specific cases $q=0$ and $\eps'_\rms=1$ which make $\ph_3$
vanish shall be treated separately. The general case is treated by first
checking that $|\ph_0(0)|<|\ph^*_0(0)|$, which is obvious. We then notice that
$\ph^*_1(0)>0$. The relation $\ph_1(0)<\ph^*_1(0)$ is equivalent to
$-\omega(\eps'_\rms-1)q<8$, which always holds. As for relation
$-\ph_1(0)<\ph^*_1(0)$, it can be cast as 
$\omega(\eps'_\rms-1)q<4\omega\eps'_\rms$, which holds true if $q\leq4$. We
therefore have $|\ph_1(0)|<|\ph^*_1(0)|$ for $q\in]0,4]$. We carry on by
studying the sign of
\begin{equation*}
\ph^*_2(0)
= \delta^2\omega[4\eps'_\rms-(\eps'_\rms-1)q][8+(\eps'_\rms-1)\omega q]>0
\end{equation*}
and therefore we have to check that $\ph^*_2(0)+\ph_2(0)>0$ and 
$\ph^*_2(0)-\ph_2(0)>0$
\begin{eqnarray*}
\ph^*_2(0)+\ph_2(0) & = & 
2\delta^2\omega[(\eps'_s-1)\omega q^2 + 32 + 2(\eps'_\rms-1)(4-q)^2] >0, \\
\ph^*_2(0)-\ph_2(0) 
& = & 2\delta^2q\omega(\eps'_\rms-1)(2+\eps'_\rms\omega)(4-q) >0,
\end{eqnarray*}
if $q\in]0,4[$. Last let us study $\ph_3$
\begin{eqnarray*}
\ph^*_3(0)
& = & 4\delta^4\omega^2(\eps'_\rms-1)q(4-q)(2+\eps'_\rms\omega)
[32+(\eps'_\rms-1)(2(4-q)+\omega q^2)]>0, \\
\ph^*_3(0)+\ph_3(0) 
& = & 8\delta^4\omega^2(\eps'_\rms-1)q^2(4-q)(2+\eps'_\rms\omega)
[8 + \omega(\eps'_s-1)q] >0, \\
\ph^*_3(0)-\ph_3(0) 
& = &  8\delta^4\omega^2(\eps'_\rms-1)q(4-q)^2(2+\eps'_\rms\omega)
[4\eps'_\rms - (\eps'_s-1)q]>0.
\end{eqnarray*}
Hence we show that $\ph_3$ and therefore $\ph_0$ is a Schur polynomial if
$q\in]0,4[$ and there remains to treat the specific cases.

\subsubsection{Case $q=0$}

The case $q=0$ has once more to be treated separately. We have
\begin{eqnarray*}
\ph_0(Z) & = & (Z-1)^2 \big[ (Z-1)^2 + \delta (Z^2-1) 
+ \frac12\omega\eps'_\rms (Z+1)^2 \big].
\end{eqnarray*} 
The corresponding amplification matrix is
\begin{equation*}
\left(\begin{array}{cccc}
1 & 0 & 0 & 0 \\
0 & \frac{D-\omega(\eps'_\rms-1)}{D} & \frac{\omega}{D} & \frac{-1}{D} \\
0 & \frac{\omega(\eps'_\rms-1)}{D} & \frac{D-\omega}{D} & \frac{1}{D} \\
0 & \frac{2\omega(\eps'_\rms-1)}{D} & \frac{-2\omega}{D} & \frac{2-D}{D}
\end{array}\right).
\end{equation*}
Anew 1 is a double eigenvalue in two distinct eigensubspaces. The other factor 
of the characteristic polynomial, namely 
\begin{equation*}
\psi_0(Z) = [1+\delta+\frac12\omega\eps'_\rms]Z^2 - [2-\omega\eps'_\rms] Z 
+ [1-\delta+\frac12\omega\eps'_\rms],
\end{equation*}
should be a Schur (or a simple von Neumann) polynomial. We do have
$|\psi_0(0)| < |\psi^*_0(0)|$ and we compute
\begin{equation*}
\psi_1(Z) = 4\delta\{[1+\frac12\omega\eps'_\rms]Z 
- [1-\frac12\omega\eps'_\rms]\}.
\end{equation*}
Both remaining eigenvalues have a strictly lower to 1 modulus and iterates of
the amplification matrix are bounded. This holds even if $\eps'_\rms=1$.

\subsubsection{Case $q=4$}

In the case when $q=4$, 
\begin{eqnarray*}
\phi_0(Z) & = & 
[1+\delta+\frac12\omega\eps'_\rms] Z^4 + [2\delta+2\omega] Z^3
+ [-2-\omega\eps'_\rms+4\omega] Z^2 + [-2\delta+2\omega] Z 
+ [1-\delta+\frac12\omega\eps'_\rms] \\
& = & 
(Z+1)^2\{[1+\delta+\frac12\omega\eps'_\rms] Z^2 
- 2[1-\omega+\frac12\omega\eps'_\rms] Z + [1-\delta+\frac12\omega\eps'_\rms]\}.
\end{eqnarray*}
We have a double root $Z=-1$. We therefore have to study the amplification
matrix which reads
\begin{equation*}
G = \left(\begin{array}{cccc}
1 & -\sigma & 0 & 0 \\
\frac{\sigma^*(D-\frac12\omega(\eps'_\rms-1))}{D} 
& \frac{-3D+\omega(\eps'_\rms-1)}{D} 
& \frac{\omega}{D} & \frac{-1}{D} \\
\frac{\sigma^*\frac12\omega(\eps'_\rms-1)}{D} 
& \frac{-2\frac12\omega(\eps'_\rms-1)}{D} 
& \frac{D-\omega}{D} & \frac{1}{D} \\
\frac{\sigma^*\omega(\eps'_\rms-1)}{D} 
& \frac{-2\omega(\eps'_\rms-1)}{D} & \frac{-2\omega}{D} 
& \frac{2-D}{D}
\end{array}\right).
\end{equation*}
Only the vector $(\sigma,2,0,0)\transp$ is an eigenvector associated to the
eigenvalue $-1$, and we have increasing iterates for $G$, whatever the study of
the other factor of the characteristic polynomial. The value $q=4$ gives rise 
to instabilities.

\subsubsection{Case $\eps'_\rms=1$}

If $\eps'_\rms=1$, the polynomial $\ph_3$ is identically zero and we shall 
study $\ph_2$ for $q\in]0,4[$
\begin{equation*}
\ph_2(Z) = 32 \delta^2 \omega \{ Z^2 - [2-q] Z + 1 \}.
\end{equation*}
This polynomial has two distinct complex conjugate roots with unit
modulus. Polynomials $\ph_2$ and therefore $\ph_0$ are both simple von Neumann
polynomials.

\subsubsection{Synthesis for the Lorentz--Kashiwa Model}

The scheme \eqref{LorentzP1} for the one-dimensional anharmonic 
Maxwell--Lorentz equations is stable with the condition
\begin{equation*}
\dt < \dx / c_\infty.
\end{equation*}


\subsection{Young Model}

\subsubsection{Model Setting}

The Maxwell system \eqref{MaxP} is closed by a discretization of 
\eqref{LorentzP}, namely 
\begin{equation}
\label{LorentzP2}
\begin{array}{rcl}
\ds \frac1\dt (P^{n+1}_j-P^n_j) & = & \ds J^{n+\frac12}, \\
\ds \frac1\dt (J^{n+\frac12}_j-J^{n-\frac12}_j)
& = & \ds - \frac\nu2 (J^{n+\frac12}_j+J^{n-\frac12}_j) 
+ \omega_1^2(\eps_\rms-\eps_\infty)\eps_0 E^n_j - \omega_1^2 P^n_j.
\end{array}
\end{equation}
The explicit version of system \eqref{MaxP}--\eqref{LorentzP2} deals once more
with the variable
\begin{equation*}
U^n_j = (c_\infty B^{n-\frac12}_{j+\frac12}, E^n_j, 
P^n_j/\eps_0\eps_\infty,\dt J^{n-\frac12}_j/\eps_0\eps_\infty)\transp
= (\calB^{n-\frac12}_{j+\frac12}, \calE^n_j, \calP^n_j,
\calJ^{n-\frac12}_j)\transp
\end{equation*}
and reads
\begin{eqnarray*}
\calB_{j+\frac12}^{n+\frac12} 
& = & \calB_{j+\frac12}^{n-\frac12} - \lambda (\calE_{j+1}^n-\calE_j^n), \\
{}[1+\delta] \calE^{n+1}_j & = & 
[1+\delta-2\omega(\eps'_\rms-1)]\calE^n_j 
+ \lambda^2 [1+\delta](\calE^n_{j+1}-2\calE_j^n+\calE_{j-1}^n)
- \lambda (\calB_{j+\frac12}^{n-\frac12} - \calB_{j-\frac12}^{n-\frac12}) \\
&& + 2\omega \calP^n_j - [1-\delta] \calJ^{n-\frac12}_j, \\
{}[1+\delta] \calP^{n+1}_j & = & [1+\delta-2\omega]\calP^n_j 
+ 2\omega(\eps'_\rms-1) \calE^n_j + [1-\delta] \calJ^{n-\frac12}_j, \\
{}[1+\delta] \calJ^{n+\frac12}_j & = & [1-\delta]\calJ^{n-\frac12}_j 
+ 2\omega(\eps'_\rms-1) \calE^n_j - 2\omega \calP^n_j,
\end{eqnarray*} 
from which stems the amplification matrix 
\begin{equation*}
G = \left(\begin{array}{cccc}
1 & -\sigma & 0 & 0 \\
\sigma^* & \frac{(1-q)(1+\delta)-2\omega\alpha}{1+\delta} 
& \frac{2\omega}{1+\delta} & -\frac{1-\delta}{1+\delta} \\
0 & \frac{2\omega\alpha}{1+\delta} 
& \frac{1+\delta-2\omega}{1+\delta} & \frac{1-\delta}{1+\delta} \\
0 & \frac{2\omega\alpha}{1+\delta} 
& \frac{-2\omega}{1+\delta} & \frac{1-\delta}{1+\delta}
\end{array}\right).
\end{equation*}

\subsubsection{Computation of the Characteristic Polynomial}

The characteristic polynomial of $G$ is
\begin{equation*}
P(Z) = \left|\begin{array}{cccc}
Z-1 & \sigma & 0 & 0 \\
-\sigma^* & Z-\frac{(1-q)(1+\delta)-2\omega\alpha}{1+\delta} 
& -\frac{2\omega}{1+\delta} & \frac{1-\delta}{1+\delta} \\
0 & -\frac{2\omega\alpha}{1+\delta} 
& Z-\frac{1+\delta-2\omega}{1+\delta} & -\frac{1-\delta}{1+\delta} \\
0 & -\frac{2\omega\alpha}{1+\delta} 
& \frac{2\omega}{1+\delta} & Z-\frac{1-\delta}{1+\delta}
\end{array}\right|. 
\end{equation*}
Therefore setting $X=(1+\delta)(Z-1)$,
\begin{eqnarray*}
(1+\delta)^4P(Z) & = & \left|\begin{array}{cccc}
X & (1+\delta)\sigma & 0 & 0 \\
-\sigma^*(1+\delta)
& X + (1+\delta)q +2\omega\alpha & -2\omega & 1-\delta \\
0 & -2\omega\alpha & X+2\omega & -(1-\delta) \\
0 & -2\omega\alpha & 2\omega & X + 2\delta
\end{array}\right| \\
& = & \left|\begin{array}{cccc}
X & (1+\delta)\sigma & 0 & 0 \\
-\sigma^*(1+\delta) & X + (1+\delta)q & X & 0 \\
0 & -2\omega\alpha & X+2\omega & -(1-\delta) \\
0 & 0 & -X & X + 1 + \delta
\end{array}\right| \allowdisplaybreaks \\
& = & X\left|\begin{array}{ccc}
X + (1+\delta)q & X & 0 \\
-2\omega\alpha & X+2\omega & -(1-\delta) \\
0 & -X & X + 1 + \delta
\end{array}\right| \\
&& + (1+\delta)^2 q \left|\begin{array}{ccc}
1 & 0 & 0 \\
-2\omega\alpha & X+2\omega & -(1-\delta) \\
0 & -X & X + 1 + \delta
\end{array}\right| \\
& = & X \{ X^3 + [2\delta+2\omega\eps'_\rms + (1+\delta)q] X^2 
+ [2(1+\delta)\omega\eps'_\rms +2(1+\delta)(\delta+\omega)q] X \\
&& \hspace{5mm} + [2(1+\delta)^2\omega q]\} 
\hspace{5mm}
+ (1+\delta)^2q\{X^2 + [2(\delta+\omega)]X + [2(1+\delta)\omega]\} \\ 
& = & X^4 + X^3 [2\delta + 2\omega\eps'_\rms + (1+\delta)q] 
+ X^2 [(1+\delta) (2\omega\eps'_\rms + (1 + 3\delta + 2\omega) q)] \\
&& + X [(1+\delta)^2 q(2\delta+4\omega)] + [2(1+\delta)^3\omega q].
\end{eqnarray*}
The characteristic polynomial is proportional to
\begin{eqnarray*}
\phi_0(Z) & = & 
[1+\delta] Z^4 - [4 + 2\delta - 2\omega\eps'_\rms - (1+\delta)q] Z^3 
+ 2 [3 - 2\omega\eps'_\rms + (\omega-1)q ] Z^2 \\
&& - [4 - 2\delta - 2\omega\eps'_\rms - (1-\delta)q] Z + [1-\delta].
\end{eqnarray*}

\subsubsection{Von Neumann Analysis}

We successively compute
\begin{eqnarray*}
\ph_1(Z) & = & 4\delta \{ Z^3 - [3 - \omega\eps'_\rms - q] Z^2
+ [3 - 2\omega\eps'_\rms + (\omega-1) q] Z - [1-\omega\eps'_\rms] \}, \\
\ph_2(Z) & = & (4\delta)^2\omega \{
[\eps'_\rms(2-\omega \eps'_\rms)] Z^2 
-[2 \eps'_\rms(2 - \omega\eps'_\rms) 
+ (\eps'_\rms(\omega-1)-1))q ] Z\\
&& \hspace{1cm} + [\eps'_\rms(2-\omega \eps'_\rms) 
- (\eps'_\rms-1)q]
\}, \\\
\ph_3(Z) & = & (4\delta)^4q\omega^2(\eps'_\rms-1) \{
[2\eps'_\rms(2-\omega\eps'_\rms) - (\eps'_\rms-1)q]Z \\
&& \hspace{3cm} 
-[2\eps'_\rms(2-\omega\eps'_\rms) - (\eps'_\rms-1)q 
- (2-\omega\eps'_\rms)q]  \}.
\end{eqnarray*}
We see that we shall once more treat the cases $q=0$ and $\eps'_\rms=1$
separately since $\ph_3$ is identically zero. Let us check the conditions in 
the general case. First, $|\phi_0(0)|<|\phi_0^*(0)|$ clearly holds as well as
$|\phi_1(0)|<|\phi_1^*(0)|$ under the condition $\omega<2/\eps'_\rms$. If 
$\eps'_s>1$ and $q\neq0$, the condition $|\phi_2(0)|<|\phi_2^*(0)|$ is 
equivalent to
\begin{equation*}
(\eps'_\rms-1)q < 2\eps'_\rms(2-\omega\eps'_\rms).
\end{equation*}
If the worst case is $q=2$, we must have
$(\eps'_\rms-1)2 < 2\eps'_\rms(2-\omega\eps'_\rms)$, which is equivalent to
$\omega < (\eps'_\rms+1)/{\eps'_\rms}^2$. If the worst case is $q=4$, we must
have $(\eps'_\rms-1)4 < 2\eps'_\rms(2-\omega\eps'_\rms)$, which is equivalent 
to $\omega < 2/{\eps'_\rms}^2$. We wait until the study of $\ph_3$ to choose
between $q\leq2$ and $q\leq4$. The root of $\phi_3$ is 
\begin{equation*}
Z = \frac{2\eps'_\rms(2-\omega\eps'_\rms) - (\eps'_\rms-1)q 
- (2-\omega\eps'_\rms)q}
{2\eps'_\rms(2-\omega\eps'_\rms) - (\eps'_\rms-1)q}.
\end{equation*}
The denominator is positive under the same assumption found to ensure 
$|\phi_2(0)|<|\phi_2^*(0)|$. If we want $|Z|<1$, the condition is hence
\begin{equation*}
(2-\omega\eps'_\rms)q < 4\eps'_\rms(2-\omega\eps'_\rms) - 2(\eps'_\rms-1)q .
\end{equation*}
If the worst case is $q=2$, we must have
$2(2-\omega\eps'_\rms) < 4\eps'_\rms(2-\omega\eps'_\rms)- 4(\eps'_\rms-1)$, 
which is equivalent to $\omega < 2/(2\eps'_\rms-1)$. If the worst case is  
$q=4$, we must have
$4(2-\omega\eps'_\rms) < 4\eps'_\rms(2-\omega\eps'_\rms) - 8(\eps'_\rms-1)$, 
which is equivalent to $4\omega\eps'_\rms(\eps'_\rms-1) < 0$, which is false. 
We therefore choose to take $q\leq2$ and the successive conditions found are
\begin{equation*}
\omega < \frac2{\eps'_\rms}, \hspace{1cm} 
\omega < \frac{\eps'_\rms+1}{{\eps'_\rms}^2}, \hspace{1cm}
\omega < \frac2{2\eps'_\rms-1}.
\end{equation*}
The more restrictive condition that we have encountered is 
$\omega< 2/(2\eps'_\rms-1)$, this is therefore our final stability condition 
in addition to $q<2$, for which $\phi_0$ is a Schur polynomial. \\

Which are the limiting case we have to study? The three conditions are
equivalent if and only if $\eps'_\rms=1$. In this case, if $\omega$ has its
limit value $\omega=2$ and $q$ its limit value $q=2$, we have 
$\ph_2(Z)\equiv0$. If $\eps'_\rms\neq1$, $q=2$ and $\omega=2/(2\eps'_\rms-1)$,
then the modulus of the root of $\ph_3$ is 1, and we conclude that $\ph_3$ and
therefore $\ph_0$ are simple von Neumann polynomials.

\subsubsection{Case $q=0$}

If $q=0$, the characteristic polynomial has the double eigenvalue 1
\begin{equation*}
\ph_0(Z) = (Z-1)^2((Z-1)^2 + \delta(Z^2-1) + 2\omega\eps'_\rms Z).
\end{equation*}
The corresponding amplification matrix is 
\begin{equation*}
G = \left(\begin{array}{cccc}
1 & 0 & 0 & 0 \\
0 & \frac{1+\delta-2\omega\alpha}{1+\delta} 
& \frac{2\omega}{1+\delta} & -\frac{1-\delta}{1+\delta} \\
0 & \frac{2\omega\alpha}{1+\delta} 
& \frac{1+\delta-2\omega}{1+\delta} & \frac{1-\delta}{1+\delta} \\
0 & \frac{2\omega\alpha}{1+\delta} 
& \frac{-2\omega}{1+\delta} & \frac{1-\delta}{1+\delta}
\end{array}\right).
\end{equation*}
The double eigenvalue operates on two distinct eigensubspaces. We have to 
study the other factor of the polynomial
\begin{equation*}
\psi_0(Z) = [1+\delta]Z^2 -[2(1-\omega\eps'_\rms)] Z + [1-\delta].
\end{equation*}
We clearly have $|\psi_0(0)|<|\psi^*_0(0)|$. Moreover we compute
\begin{equation*}
\psi_1(Z) = 4\delta\{Z-[1-\omega\eps'_\rms]\}.
\end{equation*}
We recover the condition $\omega<2/\eps'_\rms$, under which we have a Schur
polynomial. For $\omega=2/\eps'_\rms$, we have a simple von Neumann polynomial
($-1$ is a root), which allows to conclude.

\subsubsection{Case $\eps'_\rms=1$}

In the case when $\eps'_\rms=1$, only the lasts steps have to be considered. 
The only problems are the condition $|\ph_2(0)|<|\ph^*_2(0)|$ and a vanishing
$\ph_3$. In this specific case,
\begin{equation*}
\ph_2(Z) = (4\delta)^2\omega(2-\omega) \{Z^2 - [2-q] Z +1\}.
\end{equation*}
If $\omega<2/\eps'_\rms$, i.e. $\omega<2$, this polynomial is identically 
zero, and 
\begin{equation*}
\ph'_2(Z) = (4\delta)^2\omega(2-\omega) \{2Z - [2-q]\}.
\end{equation*}
Polynomial $\ph'_2$ is a Schur one if $q\in]0,2]$ and hence
$\ph_0$ is a simple von Neumann polynomial.\\
If $\omega=2$, polynomial $\ph_2$ is identically zero and we compute
\begin{equation*}
\ph_1(Z) = 4\delta (Z+1)\{Z^2 - [2-q] Z +1\}.
\end{equation*}
For $q\in]0,2]$ the roots of $Z^2 - [2-q] Z +1$ are complex conjugate, 
distinct, and their modulus is 1.

\subsubsection{Synthesis for the Lorentz--Young model}
 
The scheme \eqref{MaxP}--\eqref{LorentzP2} for the one-dimensional anharmonic
Maxwell--Lorentz is stable with the condition
\begin{equation*}
\dt \leq \min \left( \frac{\dx}{\sqrt2c_\infty}, 
\frac2{\omega_1\sqrt{2\eps'_\rms-1}}\right).
\end{equation*}

\section{Harmonic Lorentz Type Media}
\label{sec-harmonique}

Harmonic Lorentz type media are treated thanks to the three above mentioned
schemes. The computation of the amplification matrices and the characteristic
polynomials remains unchanged. The harmonicity $\nu=0$ is expressed by the
parameter $\delta=0$. This vanishing value makes $\phi_1$ identically zero for
the three schemes and the above analysis breaks down. We resume to the analysis
of the three schemes.


\subsection{Joseph et al. Model}

\subsubsection{General Case}

Since $\ph_1$ is identically zero, we want to apply 
Theorem~\ref{Th_vonNeumann}  and study the derivative polynomial of $\ph_0$, 
which we denote by
\begin{equation*}
\psi_0(Z) = [4+4\omega\eps'_\rms] Z^3 
- [12 + 6\omega\eps'_\rms - 3(1 + \omega) q] Z^2 
+ [12 + 4\omega\eps'_\rms - 4q] Z 
- [4 + 2\omega\eps'_\rms - (1 + \omega) q].
\end{equation*}
We notice that $\psi^*_0(0)>0$ and $\psi_0(0)=-(4-q)-\omega(2\eps'_\rms-q)<0$ 
for $q\leq2$. We therefore have to check that $-\psi_0(0)<\psi^*_0(0)$, which 
is equivalent to $-(\omega+1)q<2\omega\eps'_\rms$ and always holds. We 
therefore have $|\psi_0(0)|<|\psi^*_0(0)|$. Let us now compute
\begin{eqnarray*}
\psi_1(Z) & = & 
[4\omega\eps'_\rms(4+3\omega\eps'_\rms)+4(1+\omega)(2+\omega\eps'_\rms)q
-(1+\omega)^2q^2] Z^2 \\
&& + [-16\omega\eps'_\rms(2+\omega\eps'_\rms)+8(\omega^2\eps'_\rms-2)q
+4(1+\omega)q^2] Z \\
&& + [4\omega\eps'_\rms(4+\omega\eps'_\rms)
+4(2-\omega\eps'_\rms+6\omega+3\omega^2\eps'_\rms)q-3(1+\omega)^2q^2].
\end{eqnarray*}
Anew we have
\begin{equation*}
\psi^*_1(0) = 4\omega\eps'_\rms(4+3\omega\eps'_\rms)
+ (1+\omega)q[4(2+\omega\eps'_\rms)-(1+\omega)q] > 0.
\end{equation*}
Instead of studying the sign of $\psi_1(0)$, we will check that 
$-\psi_1(0)<\psi^*_1(0)$ and $\psi_1(0)<\psi^*_1(0)$. The relation
$-\psi_1(0)<\psi^*_1(0)$ is equivalent to
\begin{equation*}
4\omega\eps'_\rms(2+\omega\eps'_\rms) + (4-q)(1+\omega)^2q 
+ 4\omega^2(\eps'_\rms-1)q> 0,
\end{equation*}
which clearly holds. There remains $\psi_1(0)<\psi^*_1(0)$ which reverts to 
\begin{equation*}
4\omega^2{\eps'_\rms}^2 + 4\omega(\eps'_\rms-2-\omega\eps'_\rms) 
+ (1+\omega)^2q^2 >0. 
\end{equation*}
Cast like this it is not easy to conclude, but we can write it has a polynomial
of the variable $\omega$
\begin{equation*}
\omega^2(2\eps'_\rms-q)^2 + 2\omega q(2(\eps'_\rms-2)+q) + q^2 >0. 
\end{equation*}
The reduced discriminant of this polynomial is
\begin{equation*}
\Delta' = -8(\eps'_\rms-1)q^2(2-q) < 0, 
\end{equation*}
if $0<q<2$ and $\eps'_\rms\neq1$. The polynomial (in $\omega$) is therefore
always positive, which we were seeking. Hence we have
$|\psi_1(0)|<|\psi^*_1(0)|$ and we can profitably carry on with the computation
of $\psi_2(Z)$ which is the product
\begin{eqnarray*}
\psi_2(Z) 
& = & 8 [4\omega^2{\eps'_\rms}^2 
	+ (4\omega\eps'_\rms-4\omega^2\eps'_\rms-8\omega)q+(1+\omega)^2q^2]
\times \\
&& \times \{[4\omega^2{\eps'_\rms}^2 + 8\omega\eps'_\rms 
+ (4+4\omega\eps'_\rms+8\omega)q-(1+\omega)^2q^2]Z \\
&& \hspace{2cm} 
-[4\omega^2{\eps'_\rms}^2 + 8\omega\eps'_\rms + (4-2\omega\eps'_\rms)q
-(1+\omega)q^2]\}
\end{eqnarray*}
We notice that
\begin{equation*}
8[4\omega^2{\eps'_\rms}^2 + (4\omega\eps'_\rms-4\omega^2\eps'_\rms-8\omega)q
+(1+\omega)^2q^2] =
8[(4\eps'_\rms-q)\omega-q]^2 +16\omega q[4(\eps'_\rms-1)+2\eps'_\rms]>0
\end{equation*}
and we simplify by this factor denoting
\begin{eqnarray*}
\bar\psi_2(Z) & = & 
[4\omega^2{\eps'_\rms}^2 + 8\omega\eps'_\rms 
+ (4+4\omega^2\eps'_\rms+8\omega)q-(1+\omega)^2q^2]Z \\
&& -[4\omega^2{\eps'_\rms}^2 + 8\omega\eps'_\rms + (4-2\omega\eps'_\rms)q
-(1+\omega)q^2].
\end{eqnarray*}
We see that
\begin{equation*}
\bar\psi^*_2(0) 
= 4\omega^2[{\eps'_\rms}^2+(\eps'_\rms-1)q]+(1+\omega)^2q(4-q)\geq0
\end{equation*}
and to prove $|\psi_2(0)|<|\psi^*_2(0)|$, we only have to check that 
$\bar\psi^*_2(0)+\bar\psi_2(0)>0$ and $\bar\psi^*_2(0)-\bar\psi_2(0)>0$. We 
"notice" that
\begin{eqnarray*}
\bar\psi^*_2(0)+\bar\psi_2(0)& = & \omega q[(4\eps'_\rms-q)\omega+2(8-q)] >0,\\
\bar\psi^*_2(0)-\bar\psi_2(0)& = & [2\omega\eps'_\rms+(1+\omega)q]
[(4\eps'_s-q)\omega+2(4-q)] >0,
\end{eqnarray*}
which ends the proof in the general case for $\delta=0$.

\subsubsection{Case $q=0$}

In the case when $q=0$, along with the fact that $1$ is a double root "which
does not cause any trouble", $\ph_0$ has the same roots as those of
$[1+\omega\eps'_\rms]Z^2-2Z+[1+\omega\eps'_\rms]$ which are complex conjugate,
distinct, and their modulus is 1, if $\eps'_s\geq1$. 

\subsubsection{Case $q=2$}

The same holds for $q=2$ and this time, along with the roots $\pm i$, we have
the roots of the polynomial 
$[1+\omega\eps'_\rms]Z^2-2[1+\omega(\eps'_\rms-1)]Z+[1+\omega\eps'_\rms]$ which
has two complex conjugate, distinct roots with modulus 1, if $\eps'_s\geq1$.

\subsubsection{Case $\eps'_\rms=1$}

Finally if $\eps'_s=1$ (and $q\in]0,2[$), we shall return to the polynomial 
$\ph_0$ which can be cast as
\begin{equation*}
\ph_0(Z) = [Z^2-(2-q)Z+1][(1+\omega)Z^2-2Z+(1+\omega)].
\end{equation*}
Each of the second degree polynomials has two distinct complex conjugate
roots. We therefore have a simple von Neumann polynomial except in the
particular case when the two polynomials are proportional and have the same
roots, which are then double roots. This is reached if $(2-q) = 2/(1+\omega)$,
namely $q=2\omega/(1+\omega)$ or equivalently $\omega=q/(2-q)$. In this case,
von Neumann analysis is not useful anymore and we have to revert to the
amplification matrix
\begin{equation*}
G = \left(\begin{array}{cccc}
1 & -\sigma & 0 & 0 \\
0 & -2q+2 & -1 & q \\
0 & 1 & 0 & 0 \\
\sigma^* & -q & 0 & 1
\end{array}\right).
\end{equation*}
The two double eigenvalues of this matrix are $(2-q\pm i\sqrt{q(4-q)})/2$ and 
their each only have one corresponding eigenvector
\begin{equation*}
\left(\sigma,\frac{q\mp i\sqrt{q(4-q)}}2,
\frac{q(3-q)\mp i(2-q)\sqrt{q(4-q)}}2,
\frac{-q\pm i\sqrt{q(4-q)}}2\right)\transp.
\end{equation*}
For each eigenvalue the associated minimal eigensubspace is therefore
two-dimensional, which corresponds to an unstable case. We can say that the
scheme is stable for $q\in[0,2\omega/(1+\omega)[$. If we rewrite this in
physical variables, we have
\begin{equation*}
4c_\infty^2 \frac{\dt^2}{\dx^2}\sin^2(\frac\xi2) 
< \frac{2\omega_1^2\dt^2}{2+\omega_1^2\dt^2}.
\end{equation*}
If $\dx^2<4c_\infty^2/\omega_1^2$, this is not a bound on the time step. 
If $\dx$ if large enough, the bound on the time step is
\begin{equation*}
\dt^2 < \frac{\dx^2}{2c_\infty^2}-\frac2{\omega_1^2}.
\end{equation*}

\subsubsection{Synthesis for the Harmonic Lorentz--Joseph et al. Model}

The scheme \eqref{Max}--\eqref{LorentzD1} for the one-dimensional harmonic 
Maxwell--Lorentz equations is stable with the condition
\begin{equation*}
\dt \leq \frac{\dx}{\sqrt2 c_\infty} \textrm{ if } \eps_\rms > \eps_\infty 
\textrm{ \ \ and \ \ } 
\dt < \sqrt{\frac{\dx^2}{2c_\infty^2}-\frac2{\omega_1^2}} 
\textrm{ if } \eps_\rms = \eps_\infty,
\end{equation*}
this last condition being meaningful only if $\dx>2c_\infty/\omega_1$. It is 
therefore advisable not to use the Lorentz--Joseph et al. scheme in the 
harmonic case for $\eps_\rms=\eps_\infty$.


\subsection{Kashiwa Model}

\subsubsection{General Case}

Anew the polynomial $\ph_1$ is identically zero and we study the derivative of
polynomial $\ph_0$, which we denote
\begin{equation*}
\psi_0(Z) = [4+2\omega\eps'_\rms] Z^3 
- [12 - 3(1 + \frac12\omega) q] Z^2 
+ [12 - 2\omega\eps'_\rms - (4-2\omega)q] Z 
- [4 - (1 + \frac12\omega) q].
\end{equation*}
The condition $|\psi_0(0)|<|\psi^*_0(0)|$ is equivalent to 
$(8-q)+\frac12\omega(4\eps'_\rms-q)>0$ which holds for $q\leq4$. Then we 
compute
\begin{eqnarray*}
\psi_1(Z) & = &  
[4(4+\omega\eps'_\rms)\omega\eps'_\rms + 4(2+\omega)q
- (1+\frac12\omega)^2q^2] Z^2 \\
&& - [32\omega\eps'_\rms + (16-8\omega-8\omega\eps'_\rms-4\omega^2\eps'_\rms)q 
+ (\omega^2-4)q^2] Z \\
&& + [4(4-\omega\eps'_\rms)\omega\eps'_\rms + 4(2-2\omega\eps'_\rms+5\omega)q
- 3(1+\frac12\omega)^2q^2].
\end{eqnarray*}
We check that 
\begin{eqnarray*}
\psi^*_1(0) 
& = & \frac14 [4\omega\eps'_\rms+(2+\omega)q]
[16+4\omega\eps'_\rms-(2+\omega)q] >0,\\
\psi^*_1(0)+\psi_1(0) 
& = & q(4-q)(\omega+2)^2+4q(\eps'_\rms-1)\omega^2+8\eps'_\rms(4-q)+8q >0, \\
\psi^*_1(0)-\psi_1(0) 
& = & \frac12\{[(4\eps'_\rms-q)\omega-2q]^2+32q\omega(\eps'_\rms-1)\}>0,
\end{eqnarray*}
under the only condition that $q\in[0,4]$, which ensures that
$|\psi_1(0)|<|\psi^*_1(0)|$. Last we compute $\psi_2(Z)$ which can be cast as
\begin{equation*}
\psi_2(Z) 
= \frac12\{[(-8+4\eps'_\rms+q)\omega+2q]^2+16\alpha(4-q)\} \bar\psi_2(Z) 
\end{equation*}
with
\begin{equation*}
\bar\psi_2(Z) = 
[(2+\omega)^2q(4-q)+4((\eps'_\rms-1) \omega^2q+2(\eps'_\rms-1)\omega(4-q)+8)]Z
-[8-(2+\omega)q][4\omega\eps'_\rms+(2-\omega)q].
\end{equation*}
We check that
\begin{eqnarray*}
\bar\psi^*_2(0) & \geq & 0, \\
\bar\psi^*_2(0)+\bar\psi_2(0) & = & q\omega[(6\eps'_\rms-q)\omega+2(8-q)]>0, \\
\bar\psi^*_2(0)-\bar\psi_2(0) & = & 2(4-q)[4\eps'_\rms\omega+(2+\omega)q]>0.
\end{eqnarray*}
These two last inequalities are strict only if $q\in]0,4[$, and we then have
$|\psi_2(0)|<|\psi^*_2(0)|$. The case $\eps'_\rms=1$ is not specific in this
general study.

\subsubsection{Case $q=0$}

To treat the specific case when $q=0$, we have to revert to the study of 
$\ph_0$ which is here
\begin{eqnarray*}
\phi_0(Z) & = & 
[1+\frac12\omega\eps'_\rms] Z^4 - 4 Z^3 + [6-\omega\eps'_\rms] Z^2 - 4 Z 
+ [1+\frac12\omega\eps'_\rms] \\
& = & (Z-1)^2 \{[1+\frac12\omega\eps'_\rms] Z^2
- 2 [1-\frac12\omega\eps'_\rms]Z
+ [1+\frac12\omega\eps'_\rms]\}.
\end{eqnarray*}
The double eigenvalue $Z=1$ is the same as in the anharmonic case and does not
cause any trouble either (minimal eigensubspaces are still one-dimensional). 
The other factor of the polynomial has clearly two distinct complex conjugate 
roots of modulus 1. This configuration corresponds to a stability case for the 
scheme.

\subsubsection{Case $q=4$}

The analysis performed in the anharmonic case remains valid here. The 
eigenvalue $Z=-1$ is double and the associated minimal eigensubspace is
two-dimensional. Iterates of the amplification matrix are therefore linearly
increasing and the case is unstable.

\subsubsection{Synthesis for the Harmonic Lorentz--Kashiwa Model}

The scheme \eqref{LorentzP1} for one-dimensional harmonic Maxwell--Lorentz 
equations is stable with the condition
\begin{equation*}
\dt < \dx/c_\infty.
\end{equation*}


\subsection{Young Model}

\subsubsection{general Case}

Once more, the polynomial $\ph_1$ is identically zero and we study the 
derivative of polynomial $\ph_0$, which we denote
\begin{equation*}
\psi_0(Z) =  4 Z^3 + [-12 + 6\omega\eps'_\rms + 3q] Z^2 
+ [12 - 8\omega\eps'_\rms + 4(\omega-1)q] Z + [-4 + 2\omega\eps'_\rms + q].
\end{equation*}
The condition $|\psi_0(0)|<|\psi^*_0(0)|$ is equivalent to
$(4-2\omega\eps'_\rms)+(4-q)>0$, which we assume ($\omega<2/\eps'_\rms$ and 
$q\in[0,2]$). We carry on computing
\begin{eqnarray*}
\psi_1(Z) & = & 
[(2\omega\eps'_\rms+q)(8-(2\omega\eps'_\rms+q))]Z^2
+[(2\omega\eps'_\rms+q)(4(2\omega\eps'_\rms+q)-16-4\omega q)+16\omega q]Z \\
&& +[(2\omega\eps'_\rms+q)(8-3(2\omega\eps'_\rms+q))+16\omega q].
\end{eqnarray*}
We see immediately in this formulation that $\psi^*_1(0)>0$. Besides
\begin{eqnarray*}
\psi^*_1(0)+\psi_1(0) & = & 
4\{-4\omega^2{\eps'_\rms}^2 +8\omega\eps'_\rms+4(\omega(\eps'_\rms+1)+1)q
-q^2\},\\
\psi^*_1(0)-\psi_1(0) & = & 
2\{[2\omega-q]^2+4\omega(\eps'_\rms-1)[q+\omega(\eps'_\rms+1)\}]>0,
\end{eqnarray*}
Let us note $f(\omega)=-4\omega^2{\eps'_\rms}^2 +8\omega\eps'_\rms
+4(\omega(\eps'_\rms+1)+1)q-q^2$, we want to prove that this quantity is
positive for $\omega\in]0,2/(2\eps'_\rms-1)]$. We derive to obtain
$f'(\omega)=-4\eps'_\rms q-8\omega{\eps'_\rms}^2+8\eps'_\rms+4q$ which vanishes
at $\omega=(2\eps'_\rms-(\eps'_\rms-1)q)/2{\eps'_\rms}^2$. This corresponds to
values of $\omega$ between $1/{\eps'_\rms}^2$ (value for $q=2$) and
$1/\eps'_\rms$ (value for $q=0$), which always belong to the interval
$]0,2/(2\eps'_\rms-1)]$. At this point we have a maximum of the function
$f(\omega)$. To conclude, we only have to evaluate the limit for $f(\omega)$ as
$\omega\to0$ and the value of $f(\omega)$ at $\omega=2/(2\eps'_\rms-1)$. If 
both values are positive, $f(\omega)$ will be positive on the whole interval. 
\begin{eqnarray*}
\lim_{\omega\to0} f(\omega) & = & 4q - q^2 >0, \\
f(2/(2\eps'_\rms-1)) & = & (2\eps'_\rms-1)^2(2-q) + 2(2\eps'_\rms-1)
+ 4(\eps'_\rms-1)(2\eps'_\rms+1) > 0,
\end{eqnarray*}
If $\eps'_\rms=1$, $q=2$ and $\omega=2$, we have $\psi^*_1(0)=-\psi_1(0)$. We
will treat this case apart. We finally compute $\psi_2(Z)$ which can be cast as
\begin{equation*}
\psi_2(Z) = 8[(2\omega\eps'_\rms-q)^2+8\omega(\eps'_\rms-1)q]\bar\psi_2(Z)
\end{equation*}
with
\begin{equation*}
\bar\psi_2(Z) = [-4\omega^2{\eps'_\rms}^2 +8\omega\eps'_\rms
+4(\omega(\eps'_\rms+1)+1)q-q^2]Z
-[(2\omega\eps'_\rms+(1-\omega)q)(-4+2\omega\eps'_\rms+q)]. 
\end{equation*}
We notice that $\psi_2$ is identically zero if $\eps'_\rms=1$ and 
$q=2\omega\eps'_s=2\omega$, which has to be treated separately. In the 
opposite case, we check that
\begin{eqnarray*}
\bar\psi^*_2(0) & \geq & f(\omega) > 0, \\
\bar\psi^*_2(0)+\bar\psi_2(0) & = & q\omega[8-2\omega\eps'_\rms-q] > 0,\\
\bar\psi^*_2(0)-\bar\psi_2(0) & = & 
[2\omega\eps'_\rms+q][8-4\omega\eps'_\rms+(\omega-2)q].
\end{eqnarray*}
The quantity $8-4\omega\eps'_\rms+(\omega-2)q$ is minimum if $q=2$ (since 
$\omega<2$) and is then equal to $4-2\omega(2\eps'_\rms-1)$. As in the
anharmonic case if $\omega<2/(2\eps'_\rms-1)$ this quantity is positive, and if
$\omega=2/(2\eps'_\rms-1)$ this quantity is zero. Yet we want to show that
$\psi_0$ is a Schur or a simple von Neumann polynomial, we therefore have to
revert to the study of $\phi_0$. Once more the cases $\eps'_\rms=1$, $q=2$ and
$\omega=2$ have to be treated specifically.

\subsubsection{Case $q=0$}

In the case when $q=0$, we revert to the direct study of $\phi_0$
\begin{eqnarray*}
\phi_0(Z) & = & 
Z^4 - [4 - 2\omega\eps'_\rms] Z^3 + 2 [3 - 2\omega\eps'_\rms] Z^2 
- [4 - 2\omega\eps'_\rms] Z + 1 \\
& = & (Z-1)^2[Z^2-2(1-\omega\eps'_\rms)Z+1],
\end{eqnarray*}
which has $Z=1$ as a double root, which is no more a problem as in the
anharmonic case. Both other roots are complex conjugate, distinct and have a 
unit modulus if $\omega\eps'_\rms<2$. If $\eps'_\rms=1$ and $\omega=2$, $Z=-1$ 
is also a double root. Then we have
\begin{equation*}
G+\Id = \left(\begin{array}{cccc}
2 & -\sigma & 0 & 0 \\
\sigma^* & 2 & 1 & -\frac12 \\
0 & 0 & 1 & \frac11 \\ 
0 & 0 & -2 & 1 
\end{array}\right),
\end{equation*}
which has only $(\sigma,2,0,0)\transp$ as eigenvalue (associated to $0$). This
is a cause of instability for the scheme.

\subsubsection{Case $q=2$}

The case $q=2$ is treated by the general case except when
$\omega=2/(2\eps'_\rms-1)$. In this case $2\omega\eps'_\rms=2+\omega$, hence
\begin{eqnarray*}
\phi_0(Z) & = & Z^4 + \omega Z^3 + 2 [\omega-1] Z^2 + \omega Z + 1 \\
& = & (Z+1)^2(Z^2-(2-\omega)Z+1). 
\end{eqnarray*}
We therefore have to study the stable subspaces through the amplification 
matrix for the eigenvalue $-1$. We have 
\begin{equation*}
G+\Id = \left(\begin{array}{cccc}
2 & -\sigma & 0 & 0 \\
\sigma^* & \omega-2 & 2\omega & -1 \\
0 & 2-\omega & 2-2\omega & 1 \\ 
0 & 2-\omega & -2\omega & 2 
\end{array}\right)
\end{equation*}
which has a unique eigenvector (associated to 0), namely 
$(\sigma,2,-1,-2)\transp$.  This is an unstable case.

\subsubsection{Case $\eps'_\rms=1$}

Only the case when $\eps'_\rms=1$, $q=2\omega$ remains to study, in which case 
$\psi_2$ vanishes. We compute $\phi_0$ which is equal to
\begin{equation*}
\psi_0(Z) = (Z^2 - 2[1-\omega]Z + 1)^2. 
\end{equation*}
The two complex conjugate roots $1-\omega\pm i\sqrt{\omega(2-\omega)}$ are 
both double with modulus 1. We therefore have to study the associated stable
subspaces. The only associated eigenvectors are 
$(\sigma,\omega\mp i \sqrt{\omega(2-\omega)})\transp$ respectively and the 
associated minimal eigensubspaces are two-dimensional, which leads to 
instability.
If $\eps'_\rms=1$, we should assume $q<2\omega$, which in physical variables
reads $\dx>2c_\infty/\omega_1$ which is not a stability condition. It should
therefore be avoided to use the Lorentz--Young scheme in the harmonic case for
$\eps_\rms=\eps_\infty$ when $q$ can reach the value $2\omega$, i.e. if
$\omega\geq1$. Another way to see this condition is to give $\omega<1$ as a
stability condition if $\eps'_\rms=1$.

\subsubsection{Synthesis for the Harmonic Lorentz--Young et al. Model}

The scheme \eqref{MaxP}--\eqref{LorentzP2} for the one-dimensional harmonic 
Maxwell--Lorentz equations is stable with the condition
\begin{equation*}
\dt < \min \left( \frac{\dx}{\sqrt2c_\infty}, 
\frac2{\omega_1\sqrt{2\eps'_\rms-1}}\right) 
\textrm{ if } \eps_\rms > \eps_\infty 
\textrm{ \ \ and \ \ } 
\dt < \min \left( \frac{\dx}{\sqrt2c_\infty}, 
\frac{\sqrt2}{\omega_1}\right) 
\textrm{ if } \eps_\rms = \eps_\infty,
\end{equation*}

\section{Basic Polynomials in Dimension 1}

The previous computations lead us to define basic polynomials associated to 
each one-dimensional scheme. We will see that these polynomials will prove 
useful in higher dimensions.

\paragraph{Debye (Joseph et al.)}

\begin{eqnarray*}
P_{DJ}(Z) 
& = & [1+\delta\eps'_\rms]Z^3 - [3 + \delta\eps'_\rms - (1+\delta)q] Z^2 
+ [3-\delta\eps'_\rms - (1-\delta)q] Z - [1-\delta\eps'_\rms] \\
& = & [1+\delta\eps'_\rms]Y^3 + [2\delta\eps'_\rms + (1+\delta)q] Y^2 
+ [(1+3\delta)q] Y + [2\delta q].
\end{eqnarray*}

\paragraph{Debye (Young)}

\begin{eqnarray*}
P_{DY}(Z) & = & 
[(1+\delta\alpha)(1+\delta)] Z^3 
- [3 + \delta + \delta\alpha + 3\delta^2\alpha - (1+\delta)q] Z^2 \\
&& + [3 - \delta - \delta\alpha + 3\delta^2\alpha - (1-\delta)q] Z 
- [(1-\delta\alpha)(1-\delta)] \\
& = & [(1+\delta)(1+\delta\alpha)] Y^3 + [2\delta(1+\alpha) + (1+\delta)q] Y^2 
+ [(1+3\delta)q] Y + [2\delta q].
\end{eqnarray*}

\paragraph{Lorentz (Joseph et al.)}

\begin{eqnarray*}
P_{LJ}(Z) & = & 
[1+\delta+\omega\eps'_\rms] Z^4 
- [4 + 2\delta + 2\omega\eps'_\rms - (1 + \delta + \omega) q] Z^3 
+ [6 + 2\omega\eps'_\rms - 2q] Z^2 \\
&& - [4 - 2\delta + 2\omega\eps'_\rms - (1-\delta+\omega) q] Z 
+ [1-\delta+\omega\eps'_\rms] \\
& = & [1+\delta+\omega\eps'_\rms] Y^4 
+ [2\delta + 2\omega\eps'_\rms + (1 + \delta + \omega) q] Y^3 
+ [2\omega\eps'_\rms + (1+3\delta+3\omega)q] Y^2 \\
&& + [2 (\delta+2\omega) q] Y + [2\omega q].
\end{eqnarray*}

\paragraph{Lorentz (Kashiwa et al.)}

\begin{eqnarray*}
P_{LK}(Z) & = & 
[1+\delta+\frac12\omega\eps'_\rms] Z^4
- [4+2\delta-(1+\delta+\frac12\omega)q] Z^3
+ [6-\omega\eps'_\rms+(\omega-2)q] Z^2 \\
&&- [4-2\delta-(1-\delta+\frac12\omega)q] Z 
+ [1-\delta+\frac12\omega\eps'_\rms] \\
& = & [1+\delta+\frac12\omega\eps'_\rms] Y^4 
+ [2(\delta+\omega\eps'_\rms) + (1+\delta+\frac12\omega)q]Y^3 \\
&& + [2\omega\eps'_\rms + (1+3\delta+\frac52\omega)q]Y^2 
+ [2(\delta+2\omega)q] Y + [2\omega q].
\end{eqnarray*}

\paragraph{Lorentz (Young)}

\begin{eqnarray*}
P_{LY}(Z) & = & 
[1+\delta] Z^4 - [4 + 2\delta - 2\omega\eps'_\rms - (1+\delta)q] Z^3 
+ 2 [3 - 2\omega\eps'_\rms + (\omega-1)q ] Z^2 \\
&& - [4 - 2\delta - 2\omega\eps'_\rms - (1-\delta)q] Z + [1-\delta] \\
& = & [1+\delta] Y^4 + [2\delta + 2\omega\eps'_\rms + (1+\delta)q] Y^3 
+[2\omega\eps'_\rms + (1+3\delta+2\omega) q]  Y^2 \\
&& + [2(\delta+2\omega)q] Y + [2\omega q].
\end{eqnarray*}

\section{The Two-Dimensional Space Case}

\subsection{Maxwell Equations}

In dimension 2, the field may be decoupled into two polarisations which lead to
different schemes and also a different number of variables. We use here similar
notations as those introduced in the one-dimensional case, namely
$\lambda_x=c_\infty\dt/\dx$, $\lambda_y=c_\infty\dt/\dy$,
$\sigma_x=\lambda_x(e^{i\xi_x}-1)$, $\sigma_y=\lambda_x(e^{i\xi_y}-1)$,
$q_x=|\sigma_x|^2$ and $q_y=|\sigma_y|^2$.

\subsubsection{Polarisation $(B_x,B_y,E_z)$}
The polarisation $(B_x,B_y,E_z)$ is also called the transverse electric 
polarisation $TE_z$.
\begin{eqnarray*}
\calB_{x,j,k+\frac12}^{n+\frac12} - \calB_{x,j,k+\frac12}^{n-\frac12}
& = & -\ly \left( \calE_{z,j,k+1}^n - \calE_{z,j,k}^n \right), \\
\calB_{y,j+\frac12,k}^{n+\frac12} - \calB_{y,j+\frac12,k}^{n-\frac12}
& = & \lx \left( \calE_{z,j+1,k}^n - \calE_{z,j,k}^n \right), \\
\calD_{z,j,k}^{n+1}-\calD_{z,j,k}^n & = & 
\lx \left( \calB_{y,j+\frac12,k}^{n+\frac12} 
- \calB_{y,j-\frac12,k}^{n+\frac12} \right)
- \ly \left( \calB_{x,j,k+\frac12}^{n+\frac12}
- \calB_{x,j,k-\frac12}^{n+\frac12} \right).
\end{eqnarray*}

\subsubsection{Polarisation $(B_z,E_x,E_y)$}

The polarisation $(B_z,E_x,E_y)$ is also called the transverse magnetic 
polarisation $TM_z$.
\begin{eqnarray*}
\calB_{z,j+\frac12,k+\frac12}^{n+\frac12}
- \calB_{z,j+\frac12,k+\frac12}^{n-\frac12} & = & 
- \lx \left( \calE_{y,j+1,k+\frac12}^n - \calE_{y,j,k+\frac12}^n \right)
+ \ly \left( \calE_{x,j+\frac12,k+1}^n - \calE_{x,j+\frac12,k}^n \right), \\
\calD_{x,j+\frac12,k}^{n+1} - \calD_{x,j+\frac12,k}^n 
& = & \ly \left( \calB_{z,j+\frac12,k+\frac12}^{n+\frac12}
- \calB_{z,j+\frac12,k-\frac12}^{n+\frac12} \right), \\
\calD_{y,j,k+\frac12}^{n+1} - \calD_{y,j,k+\frac12}^n 
& = & - \lx \left( \calB_{z,j+\frac12,k+\frac12}^{n+\frac12}
- \calB_{z,j-\frac12,k+\frac12}^{n+\frac12} \right).
\end{eqnarray*}

\subsection{The Debye--Joseph et al. Scheme}

\subsubsection{Polarisation $(B_x,B_y,E_z)$}

Coupling polarisation $(B_x,B_y,E_z)$ with the Debye--Joseph et al. scheme, we 
obtain the amplification matrix
\begin{equation*}
G = \left(\begin{array}{cccc}
1 & 0 & -\sigma_y & 0 \\
0 & 1 & \sigma_x & 0 \\
\frac{(1+\delta)\sigma_y^*}{1+\delta\eps'_\rms} 
& -\frac{(1+\delta)\sigma_x^*}{1+\delta\eps'_\rms} 
& \frac{(1-\delta\eps'_\rms)-(1+\delta)(q_x+q_y)}{1+\delta\eps'_\rms} 
& \frac{2\delta}{1+\delta\eps'_\rms} \\
\sigma_y^* & -\sigma_x^* & -(q_x+q_y) & 1 
\end{array}\right)
\end{equation*}
associated to the variable
$(\calB_{x,j,k+\frac12}^{n-\frac12},\calB_{y,j+\frac12,k}^{n-\frac12},
\calE_{z,j,k}^n,\calD_{z,j,k}^n)\transp$.
The characteristic polynomial is proportional to the characteristic polynomial
in dimension 1 for the same scheme with 1 as an extra root
\begin{equation*}
\phi_0(Z) = Y P_{DJ}(Z).
\end{equation*}
In polynomial $P_{DJ}(Z)$, the variable $q$ means $q=q_x+q_y$ in the
two-dimensional case. The polynomial only depends on this sum and not on the
separate values of $q_x$ and $q_y$. \\

The general case treated in dimension one concludes to a Schur polynomial, 
we therefore have a von Neumann polynomial here. We also see easily on the
amplification matrix in the case $q=0$ (if $q_x$ only or $q_y$ only vanish, we
do not have a specific case), that the eigensubspaces associated to the
eigenvalue 1 are indeed stable. As for the particular case $\eps'_\rms=1$, 
which gave rise to two complex conjugate eigenvalues, different from 1, we may
add this new eigenvalue with the same conclusion, namely stability with the
condition $q=q_x+q_y<4$ i.e. $\sqrt2 c_\infty \dt < \dx$ if $\dx=\dy$.

\subsubsection{Polarisation $(B_z,E_x,E_y)$}

Coupling polarisation $(B_z,E_x,E_y)$ with the Debye--Joseph et al. scheme 
yields the amplification matrix
\begin{equation*}
G = \left(\begin{array}{ccccc}
1 & \sigma_y & 0 & -\sigma_x & 0 \\
-\frac{(1+\delta)\sigma_y^*}{1+\delta\eps'_\rms} 
& \frac{(1-\delta\eps'_\rms)-(1+\delta)q_y}{1+\delta\eps'_\rms} 
& \frac{2\delta}{1+\delta\eps'_\rms} 
& \frac{(1+\delta)\sigma_x\sigma_y^*}{1+\delta\eps'_\rms} & 0 \\
- \sigma_y^* & -q_y & 1 & \sigma_x\sigma_y^* & 0 \\
\frac{(1+\delta)\sigma_x^*}{1+\delta\eps'_\rms} 
& \frac{(1+\delta)\sigma_x^*\sigma_y}{1+\delta\eps'_\rms} 
& 0 & \frac{(1-\delta\eps'_\rms)-(1+\delta)q_x}{1+\delta\eps'_\rms} 
& \frac{2\delta}{1+\delta\eps'_\rms}\\
\sigma_x^* & \sigma_x^*\sigma_y & 0 & -q_x & 1 
\end{array}\right)
\end{equation*}
associated to the variable
$(\calB_{z,j+\frac12,k+\frac12}^{n-\frac12},\calE_{x,j+\frac12,k}^n,
\calD_{x,j+\frac12,k}^n,\calE_{y,j,k+\frac12}^n,
\calD_{y,j,k+\frac12}^n)\transp$.
Anew we have a proportional polynomial to that of the one-dimensional case with
two extra roots
\begin{equation*}
\phi_0(Z) = Y[(1+\delta\eps'_\rms)Y+2\delta\eps'_\rms] P_{DJ}(Z),
\end{equation*}
which are equal to 1 and $(1-\delta\eps'_\rms)/(1+\delta\eps'_\rms)$
respectively and which are not roots in dimension 1 in the general case. \\

Only the root 1 could induce stability problems and we have seen that only the
case $q=0$ would make it a multiple root. The matrix $G$ is then 
block-diagonal, with two rank-two blocks identical to that of the 
one-dimensional case. The stability is once more ensured with the condition 
$q=q_x+q_y<4$.

\subsection{The Debye--Young Scheme}

\subsubsection{Polarisation $(B_x,B_y,E_z)$}

Coupling polarisation $(B_x,B_y,E_z)$ with Debye--Young scheme, we obtain the 
amplification matrix
\begin{equation*}
G = \left(\begin{array}{cccc}
1 & 0 & -\sigma_y & 0 \\
0 & 1 & \sigma_x & 0 \\
\frac{\sigma_y^*}{1+\delta\alpha} 
& -\frac{\sigma_x^*}{1+\delta\alpha} 
& \frac{(1+\delta)(1-\delta\alpha)+ 4\delta^2\alpha-(1+\delta)q}
{(1+\delta)(1+\delta\alpha)} 
& \frac{1-\delta}{1+\delta}\frac{2\delta}{1+\delta\alpha} \\
0 & 0 & \frac{2\delta\alpha}{1+\delta}  & \frac{1-\delta}{1+\delta} 
\end{array}\right)
\end{equation*}
associated to the variable 
$(\calB_{x,j,k+\frac12}^{n-\frac12},\calB_{y,j+\frac12,k}^{n-\frac12},
\calE_{z,j,k}^n,\calP_{z,j,k}^n)\transp$.
The computation of the characteristic polynomial leads to the same polynomial 
as in one dimension but with 1 as an extra eigenvalue
\begin{equation*}
\phi_0(Z) = Y P_{DY}(Z).
\end{equation*}
Anew we revert to the arguments of the one-dimensional case, the double 
eigenvalue 1 of the $q=0$ case not being a problem.

\subsubsection{Polarisation $(B_z,E_x,E_y)$}

Coupling polarisation $(B_z,E_x,E_y)$ with the Debye--Young scheme, we obtain 
the amplification matrix
\begin{equation*}
G = \left(\begin{array}{ccccc}
1 & \sigma_y & 0 & -\sigma_x & 0 \\
-\frac{\sigma_y^*}{1+\delta\alpha} 
& \frac{(1+\delta)(1-\delta\alpha)+4\delta^2\alpha - (1+\delta) q_y}
{(1+\delta)(1+\delta\alpha)}
& \frac{1-\delta}{1+\delta}\frac{2\delta}{1+\delta\alpha} 
& \frac{\sigma_x\sigma_y^*}{1+\delta\alpha} & 0 \\
0 & \frac{2\delta\alpha}{1+\delta} & \frac{1-\delta}{1+\delta} & 0 & 0 \\
\frac{\sigma_x^*}{1+\delta\alpha} 
& \frac{\sigma_x^*\sigma_y}{1+\delta\alpha} & 0 
& \frac{(1+\delta)(1-\delta\alpha)+4\delta^2\alpha - (1+\delta) q_x}
{(1+\delta)(1+\delta\alpha)}
& \frac{1-\delta}{1+\delta}\frac{2\delta}{1+\delta\alpha} \\
0 & 0 & 0 & \frac{2\delta\alpha}{1+\delta} & \frac{1-\delta}{1+\delta}
\end{array}\right)
\end{equation*}
associated to the variable 
$(\calB_{z,j+\frac12,k+\frac12}^{n-\frac12},\calE_{x,j+\frac12,k}^n,
\calP_{x,j+\frac12,k}^n,\calE_{y,j,k+\frac12}^n,
\calP_{y,j,k+\frac12}^n)\transp$.
The computation of the characteristic polynomial leads to the same polynomial 
as in dimension 1 but with two extra eigenvalues
\begin{equation*}
\phi_0(Z) = Y[(1+\delta)(1+\delta\alpha)Y+2\delta(1+\alpha)] P_{DY}(Z)
\end{equation*}
which are 1 and $(1-\delta)(1-\delta\alpha)/(1+\delta)(1+\delta\alpha)$.
Anew we revert to the argument in the one-dimensional case, the triple
eigenvalue 1 of case $q=0$ not being a problem.

\subsection{The Lorentz--Joseph et al. Scheme}

\subsubsection{Polarisation $(B_x,B_y,E_z)$}

Coupling polarisation $(B_x,B_y,E_z)$ with the Lorentz--Joseph et al. scheme, 
we obtain the amplification matrix
\begin{equation*}
G = \left(\begin{array}{ccccc}
1 & 0 & -\sigma_y & 0 & 0 \\
0 & 1 & \sigma_x & 0 & 0 \\
\frac{2\delta\sigma_y^*}{1+\delta+\omega\eps'_\rms} 
& -\frac{2\delta\sigma_x^*}{1+\delta+\omega\eps'_\rms} 
& \frac{2 - (1+\delta+\omega)q}{(1+\delta+\omega\eps'_\rms)} 
& -\frac{1-\delta+\omega\eps'_\rms}{1+\delta+\omega\eps'_\rms} 
& \frac{2\omega}{1+\delta+\omega\eps'_\rms}\\
0 & 0 & 1 & 0 & 0 \\
\sigma_y^* & - \sigma_x^* & -q & 0 & 1 
\end{array}\right)
\end{equation*}
associated to the variable 
$(\calB_{x,j,k+\frac12}^{n-\frac12},\calB_{y,j+\frac12,k}^{n-\frac12},
\calE_{z,j,k}^n,\calE_{z,j,k}^{n-1},\calD_{z,j,k}^n)\transp$.
The computation of the characteristic polynomial leads to the same polynomial 
as in dimension 1 but with the extra eigenvalue 1
\begin{equation*}
\phi_0(Z) = Y P_{LJ}(Z).
\end{equation*}
The triple eigenvalue in the $q=0$ case is not a problem.

\subsubsection{Polarisation $(B_z,E_x,E_y)$}

Coupling polarisation $(B_z,E_x,E_y)$ with the Lorentz--Joseph et al. scheme, 
we obtain the amplification matrix
\begin{equation*}
G = \left(\begin{array}{ccccccc}
1 & \sigma_y & 0 & 0 & -\sigma_x & 0 & 0 \\
-\frac{2\delta\sigma_y^*}{1+\delta+\omega\eps'_\rms} 
& \frac{2-(1+\delta+\omega) q_y}{1+\delta+\omega\eps'_\rms} 
& -\frac{1-\delta+\omega\eps'_\rms}{1+\delta+\omega\eps'_\rms}
& \frac{2\omega}{1+\delta+\omega\eps'_\rms}  
& \frac{\sigma_x\sigma_y^*(1+\delta+\omega)}{1+\delta+\omega\eps'_\rms} 
& 0 & 0 \\
0 & 1 & 0 & 0 & 0 & 0 & 0 \\
-\sigma_y^* & - q_y & 0 & 1 & \sigma_x\sigma_y^* & 0 & 0 \\
\frac{2\delta\sigma_x^*}{1+\delta+\omega\eps'_\rms} 
& \frac{\sigma_x^*\sigma_y(1+\delta+\omega)}{1+\delta+\omega\eps'_\rms} & 0 & 0
& \frac{2-(1+\delta+\omega) q_x}{1+\delta+\omega\eps'_\rms}
& -\frac{1-\delta+\omega\eps'_\rms}{1+\delta+\omega\eps'_\rms}
& \frac{2\omega}{1+\delta+\omega\eps'_\rms}\\
0 & 0 & 0 & 0 & 1 & 0 & 0 \\
\sigma_x^* & \sigma_x^*\sigma_y & 0 & 0 & -q_x & 0 & 1
\end{array}\right)
\end{equation*}
associated to the variable 
$(\calB_{z,j+\frac12,k+\frac12}^{n-\frac12},\calE_{x,j+\frac12,k}^n,
\calE_{x,j+\frac12,k}^{n-1},\calD_{x,j+\frac12,k}^n,\calE_{y,j,k+\frac12}^n,
\calE_{y,j,k+\frac12}^{n-1},\calD_{y,j,k+\frac12}^n)\transp$.
The characteristic polynomial is once more proportional to the one-dimensional 
polynomial
\begin{eqnarray*}
\phi_0(Z) & = & Y 
[(1+\delta+\omega\eps'_\rms)Y^2 + 2(\delta+\omega\eps'_\rms) Y 
+ 2\omega\eps'_\rms] P_{LJ}(Z)\\
& = & Y [(1+\delta+\omega\eps'_\rms)Z^2 - 2 Z + (1-\delta+\omega\eps'_\rms)] 
P_{LJ}(Z).
\end{eqnarray*}

In the anharmonic case, and by the von Neumann technique, we check easily that
\begin{equation*}
\psi_0(Z) = [1+\delta+\omega\eps'_\rms]Z^2 - 2 Z + [1-\delta+\omega\eps'_\rms] 
\end{equation*}
is a Schur polynomial. Besides the double root 1 if $q=0$ is still no 
problem. \\

In the harmonic case, $\psi_0(Z)$ has two distinct complex conjugate roots with
modulus 1. This is not a problem in itself, except if $\eps'_\rms=1$ in which 
case $\psi_0(Z)$ is also a factor in $P_{LJ}(Z)$
\begin{eqnarray*}
\phi_0(Z) & = & 
(Z-1) 
[(1+\omega\eps'_\rms)Z^2 - 2 Z + (1+\omega\eps'_\rms)] 
P_{LJ}(Z) \\
& = & (Z-1) 
[(1+\omega\eps'_\rms)Z^2 - 2 Z + (1+\omega\eps'_\rms)]^2 
[Z^2-(2-q)Z+1].
\end{eqnarray*}
This is already the case when we have detected double eigenvalues in dimension
1, giving rise to instabilities for $q=2\omega/(1+\omega)$. It is better to
avoid this scheme in the case when $\eps'_\rms$. 

\subsection{The Lorentz--Kashiwa et al. Scheme}

\subsubsection{Polarisation $(B_x,B_y,E_z)$}

Coupling polarisation $(B_x,B_y,E_z)$ with the Lorentz--Kashiwa et al. scheme, 
we obtain the amplification matrix
\begin{equation*}
G = \left(\begin{array}{ccccc}
1 & 0 & -\sigma_y & 0 & 0 \\
0 & 1 & \sigma_x & 0 & 0 \\
\frac{\sigma_y^*(D-\frac12\omega\alpha)}{D} 
& -\frac{\sigma_x^*(D-\frac12\omega\alpha)}{D} 
& \frac{(1-q)D-(2-q)\frac12\omega\alpha}{D} 
& \frac{\omega}{D} 
& -\frac{1}{D} \\
\frac{\sigma_y^*\frac12\omega\alpha}{D} 
& -\frac{\sigma_x^*\frac12\omega\alpha}{D} 
& \frac{(2-q)\frac12\omega\alpha}{D} 
& \frac{D-\omega}{D} 
& \frac{1}{D} \\
\frac{\sigma_y^*\omega\alpha}{D} 
& - \frac{\sigma_x^*\omega\alpha}{D} & 
\frac{(2-q)\omega\alpha}{D} & -\frac{2\omega}{D} & \frac{2-D}{D} 
\end{array}\right)
\end{equation*}
associated to the variable 
$(\calB_{x,j,k+\frac12}^{n-\frac12},\calB_{y,j+\frac12,k}^{n-\frac12},
\calE_{z,j,k}^n,\calP_{z,j,k}^{n},\calJ_{z,j,k}^n)\transp$.
The calculation of the characteristic polynomial leads to the same polynomial 
as in dimension 1 with the extra root 1
\begin{equation*}
\phi_0(Z) = Y P_{LK}(Z).
\end{equation*}
The triple eigenvalue 1 of case $q=0$ is not a problem.

\subsubsection{Polarisation $(B_z,E_x,E_y)$}

Coupling polarisation $(B_z,E_x,E_y)$ with the Lorentz--Kashiwa et al. scheme, 
we obtain the amplification matrix
\begin{equation*}
G = \left(\begin{array}{ccccccc}
1 & -\sigma_y & 0 & 0 & \sigma_x & 0 & 0\\
\frac{\sigma_y^*(D-\frac12\omega\alpha)}{D} 
& \frac{(1-q_y)D-(2-q_y)\frac12\omega\alpha}{D} 
& \frac{\omega}{D} 
& -\frac{1}{D}
& \frac{\sigma_x\sigma_y^*(D-\frac12\omega\alpha)}{D}
& 0 & 0 \\ 
\frac{\sigma_y^*\frac12\omega\alpha}{D} 
& \frac{(1-q_y)\frac12\omega\alpha}{D} 
& \frac{D-\omega}{D} 
& \frac{1}{D}
& \frac{\sigma_x\sigma_y^*\frac12\omega\alpha}{D}
& 0 & 0 \\ 
\frac{\sigma_y^*\omega\alpha}{D} 
& \frac{(1-q_y)\omega\alpha}{D} 
& \frac{-2\omega}{D} 
& \frac{2-D}{D}
& \frac{\sigma_x\sigma_y^*\omega\alpha}{D}
& 0 & 0 \\ 
-\frac{\sigma_x^*(D-\frac12\omega\alpha)}{D} 
& \frac{\sigma_x^*\sigma_y(D-\frac12\omega\alpha)}{D} 
& 0 & 0
& \frac{(1-q_x)D-(2-q_x)\frac12\omega\alpha}{D} 
& \frac{\omega}{D} 
& -\frac{1}{D} \\
-\frac{\sigma_x^*\frac12\omega\alpha}{D} 
& \frac{\sigma_x^*\sigma_y\frac12\omega\alpha}{D} 
& 0 & 0
& \frac{(2-q_x)\frac12\omega\alpha}{D} 
& \frac{D-\omega}{D} 
& \frac{1}{D} \\
-\frac{\sigma_x^*\omega\alpha}{D} 
& \frac{\sigma_x^*\sigma_y\omega\alpha}{D} 
& 0 & 0
& \frac{(2-q_x)\omega\alpha}{D} 
& \frac{-2\omega}{D} 
& \frac{2-D}{D} \\
\end{array}\right)
\end{equation*}
associated to the variable 
$(\calB_{z,j+\frac12,k+\frac12}^{n-\frac12},\calE_{x,j+\frac12,k}^n,
\calP_{x,j+\frac12,k}^{n},\calJ_{x,j+\frac12,k}^n,\calE_{y,j,k+\frac12}^n,
\calP_{y,j,k+\frac12}^{n},\calJ_{y,j,k+\frac12}^n)\transp$.
The computation of the characteristic polynomial leads to a polynomial
proportional to the one-dimensional one
\begin{eqnarray*}
\phi_0(Z) & = & Y [(1+\delta+\frac12\omega\eps'_\rms)Y^2 
+ (2(\delta+\omega\eps'_\rms))Y + (2\omega\eps'_\rms)]P_{LK}(Z) \\
& = & Y [(1+\delta+\frac12\omega\eps'_\rms)Z^2 
- (2-\omega\eps'_\rms)Z + (1-\delta+\frac12\omega\eps'_\rms)]P_{LK}(Z).
\end{eqnarray*}

In the anharmonic case, and by the von Neumann technique, we check easily that
\begin{equation*}
\psi_0(Z) = [1+\delta+\frac12\omega\eps'_\rms]Z^2 
- [2-\omega\eps'_\rms]Z + [1-\delta+\frac12\omega\eps'_\rms]
\end{equation*}
is a Schur polynomial. Besides, the root 1, which is a double one if $q=0$, 
does not lead to any problem. \\

In the harmonic case, we have the extra roots 1 and two complex conjugate roots
of modulus 1, which are not roots of $P_{LK}(Z)$. The stability is hence given
under the same conditions as in the one-dimensional case.

\subsection{The Lorentz--Young Scheme}

\subsubsection{Polarisation $(B_x,B_y,E_z)$}

Coupling polarisation $(B_x,B_y,E_z)$ with the Lorentz--Young scheme, we 
obtain the amplification matrix
\begin{equation*}
G = \left(\begin{array}{ccccc}
1 & 0 & -\sigma_y & 0 & 0 \\
0 & 1 & \sigma_x & 0 & 0 \\
\sigma_y^* & -\sigma_x^*
& \frac{(1+\delta)(1-q)-2\omega\alpha}{1+\delta} 
& \frac{2\omega}{1+\delta} 
& -\frac{1-\delta}{1+\delta} \\
0 & 0 & \frac{2\omega\alpha}{1+\delta} 
& \frac{1+\delta-2\omega}{1+\delta} & \frac{1-\delta}{1+\delta} \\
0 & 0 & \frac{2\omega\alpha}{1+\delta} 
& \frac{-2\omega}{1+\delta} & \frac{1-\delta}{1+\delta} \\
\end{array}\right)
\end{equation*}
associated to the variable 
$(\calB_{x,j,k+\frac12}^{n-\frac12},\calB_{y,j+\frac12,k}^{n-\frac12},
\calE_{z,j,k}^n,\calP_{z,j,k}^{n},\calJ_{z,j,k}^{n-\frac12})\transp$.
The computation of the characteristic polynomial leads to the same polynomial 
as in dimension 1 but with 1 as an extra eigenvalue
\begin{equation*}
\phi_0(Z) = Y P_{LY}(Z).
\end{equation*}
The triple eigenvalue 1 of the case $q=0$ is not a problem.

\subsubsection{Polarisation $(B_z,E_x,E_y)$}

Coupling polarisation $(B_z,E_x,E_y)$ with the Lorentz--Young scheme, we 
obtain the amplification matrix
\begin{equation*}
G = \left(\begin{array}{ccccccc}
1 & \sigma_y & 0 & 0 & - \sigma_x & 0 & 0 \\
- \sigma_y^*& \frac{(1+\delta)(1-q_y)-2\omega\alpha}{1+\delta} 
& \frac{2\omega}{1+\delta} 
& -\frac{1-\delta}{1+\delta} & \sigma_x\sigma_y^* & 0 & 0 \\
0 & \frac{2\omega\alpha}{1+\delta} & \frac{1+\delta-2\omega}{1+\delta} 
& \frac{1-\delta}{1+\delta} & 0 & 0 & 0 \\
0 & \frac{2\omega\alpha}{1+\delta} & \frac{-2\omega}{1+\delta} 
& \frac{1-\delta}{1+\delta} & 0 & 0 & 0 \\
\sigma_x^* & \sigma_x^*\sigma_y & 0 & 0 
& \frac{(1+\delta)(1-q_x)-2\omega\alpha}{1+\delta} 
& \frac{2\omega}{1+\delta} & -\frac{1-\delta}{1+\delta} \\
0 & 0 & 0 & 0 & \frac{2\omega\alpha}{1+\delta} 
& \frac{1+\delta-2\omega}{1+\delta} & \frac{1-\delta}{1+\delta} \\
0 & 0 & 0 & 0 & \frac{2\omega\alpha}{1+\delta} & \frac{-2\omega}{1+\delta} 
& \frac{1-\delta}{1+\delta} 
\end{array}\right)
\end{equation*}
associated to the variable 
$(\calB_{z,j+\frac12,k+\frac12}^{n-\frac12},
\calE_{x,j+\frac12,k}^n,\calP_{x,j+\frac12,k}^{n},
\calJ_{x,j+\frac12,k}^{n-\frac12}, \calE_{y,j,k+\frac12}^n,
\calP_{y,j,k+\frac12}^{n},\calJ_{y,j,k+\frac12}^{n-\frac12})\transp$.
The computation of the characteristic polynomial leads to a polynomial
proportional to that of dimension 1
\begin{eqnarray*}
\phi_0(Z) & = & Y [(1+\delta)Y^2 
+ (2(\delta+\omega\eps'_\rms))Y + (2\omega\eps'_\rms)]P_{LY}(Z) \\
& = & Y [(1+\delta)Z^2 - (2-\omega\eps'_\rms)Z + (1-\delta)]P_{LY}(Z).
\end{eqnarray*}

In the anharmonic case, and by von Neumann technique, we check easily that 
\begin{equation*}
\psi_0(Z) = [1+\delta]Z^2 - [2-2\omega\eps'_\rms]Z + [1-\delta]
\end{equation*}
is a Schur polynomial. Besides, the root 1 which is a double one if $q=0$ is 
not a problem. \\

In the harmonic case, we have the extra roots 1 and two complex conjugate roots
of modulus 1, which are not roots of $P_{LY}(Z)$. The stability is therefore
ensured under the same conditions as in the one-dimensional case.

\section{Conclusion}

We have studied the stability of numerical schemes for Maxwell--Debye and
Maxwell--Lorentz equations in space dimension 1 and 2. In dimension 2, the
characteristic polynomials of each scheme and in both polarisation happen to be
proportional to the characteristic polynomials for the same scheme in space
dimension 1. In all the cases, the extension to dimension 2 goes with an extra
root 1 compared to the one-dimensional case. This is the only extra root in the
$TE_z$ polarisation. For the $TM_z$ polarisation, there is one other extra root
for the Debye equation and two other extra roots for the Lorentz equation, all
these roots being on the unit circle. For the Yee scheme applied to the raw
Maxwell equations, the stability condition is $q\leq4$ in dimensions 1, 2 et 3,
recalling that $q=q_x+q_y$ in dimension 2 ($q=\max(q_x+q_y,q_x+q_z,q_y+q_z)$ in
dimension 3). The results are gathered in two tables according to
$\eps_\rms=\eps_\infty$ or not.

\begin{table}[h]
\renewcommand{\arraystretch}{1.5}
\[
\begin{array}{|c|c|c|c|c|}
\hline
\textrm{Model} & \textrm{Scheme} & 
& \textrm{dimension 1} 
& \textrm{dimension 2 ($\dx=\dy$)}\\
\hline\hline
\textrm{Debye} & \textrm{Joseph et al.} & q\leq4 
& \dt\leq\frac{\dx}{c_\infty} 
& \dt\leq\frac{\dx}{\sqrt2c_\infty} \\
\hline
\textrm{Debye} & \textrm{Young} & q\leq4,\ \delta\leq1 
& \dt\leq\min(\frac{\dx}{c_\infty},2t_\rmr) 
& \dt\leq\min(\frac{\dx}{\sqrt2c_\infty},2t_\rmr)\\
\hline\hline
\textrm{Lorentz} & \textrm{Joseph et al.} & q\leq2 
& \dt\leq\frac{\dx}{\sqrt2 c_\infty} 
& \dt\leq\frac{\dx}{2 c_\infty}\\
\hline
\textrm{Lorentz} & \textrm{Kashiwa et al.} & q<4 
& \dt<\frac{\dx}{c_\infty} 
& \dt<\frac{\dx}{\sqrt2c_\infty}\\
\hline
\textrm{Lorentz} & \textrm{Young} & q\leq2,\ \omega\leq\frac2{2\eps'_\rms-1} 
& \dt\leq\min(\frac{\dx}{\sqrt2c_\infty},\frac2{\omega_1\sqrt{2\eps'_\rms-1}}) 
& \dt\leq\min(\frac{\dx}{2 c_\infty},\frac2{\omega_1\sqrt{2\eps'_\rms-1}})\\
\hline\hline
\textrm{Harm.} & \textrm{Joseph et al.} & q\leq2 
& \dt\leq\frac{\dx}{\sqrt2 c_\infty} 
& \dt\leq\frac{\dx}{2 c_\infty}\\
\hline
\textrm{Harm.} & \textrm{Kashiwa et al.} & q<4 
& \dt<\frac{\dx}{c_\infty} 
& \dt<\frac{\dx}{\sqrt2c_\infty}\\
\hline
\textrm{Harm.} & \textrm{Young} & 
\begin{array}{c}q<2,\ \omega\leq\frac2{2\eps'_\rms-1} \\[-3mm]
\textrm{or}\\[-3mm]
q\leq2,\ \omega<\frac2{2\eps'_\rms-1} \end{array}
& \dt<\min(\frac{\dx}{\sqrt2 c_\infty},\frac2{\omega_1\sqrt{2\eps'_\rms-1}}) 
& \dt<\min(\frac{\dx}{2 c_\infty},\frac2{\omega_1\sqrt{2\eps'_\rms-1}})\\
\hline
\end{array}
\]
\caption{Stability of schemes for $\eps_\rms>\eps_\infty$.}
\end{table}

\begin{table}[h]
\renewcommand{\arraystretch}{1.5}
\[
\begin{array}{|c|c|c|c|c|}
\hline
\textrm{Model} & \textrm{Scheme} & 
& \textrm{dimension 1} 
& \textrm{dimension 2 ($\dx=\dy$)}\\
\hline\hline
\textrm{Debye} & \textrm{Joseph et al.} & q<4 
& \dt<\frac{\dx}{c_\infty} 
& \dt<\frac{\dx}{\sqrt2c_\infty} \\
\hline
\textrm{Debye} & \textrm{Young} & q<4
& \dt<\frac{\dx}{c_\infty}
& \dt<\frac{\dx}{\sqrt2c_\infty} \\
\hline\hline
\textrm{Lorentz} & \textrm{Joseph et al.} & q\leq2 
& \dt\leq\frac{\dx}{\sqrt2 c_\infty} 
& \dt\leq\frac{\dx}{2 c_\infty}\\
\hline
\textrm{Lorentz} & \textrm{Kashiwa et al.} & q<4 
& \dt<\frac{\dx}{c_\infty} 
& \dt<\frac{\dx}{\sqrt2c_\infty}\\
\hline
\textrm{Lorentz} & \textrm{Young} & q\leq2,\ \omega\leq2 
& \dt\leq\min(\frac{\dx}{\sqrt2 c_\infty},\frac2{\omega_1}) 
& \dt\leq\min(\frac{\dx}{2 c_\infty},\frac2{\omega_1})\\
\hline\hline
\textrm{Harm.} & \textrm{Joseph et al.} & q<\frac{2\omega}{1+\omega} 
& \textrm{to avoid}
& \textrm{to avoid}\\
\hline
\textrm{Harm.} & \textrm{Kashiwa et al.} & q<4 
& \dt<\frac{\dx}{c_\infty} 
& \dt<\frac{\dx}{\sqrt2c_\infty}\\
\hline
\textrm{Harm.} & \textrm{Young} & q<2,\ \omega<1
& \dt<\min(\frac{\dx}{\sqrt2 c_\infty},\frac{\sqrt2}{\omega_1}) 
& \dt<\min(\frac{\dx}{2 c_\infty},\frac{\sqrt2}{\omega_1})\\
\hline
\end{array}
\]
\caption{Stability of schemes for $\eps_\rms=\eps_\infty$.}
\end{table}

For each model, we have at least one scheme for which the stability condition 
is the same as for the raw Maxwell equations ($q<4$). In Young models, the 
extra conditions correspond to a fine enough discretization of Debye and 
Lorentz equations respectively,... because stability is not the only issue. 
Applications to classical materials show in general that the condition due to 
the Maxwell equations is the more restrictive one and not conditions due to 
the constitutive law of the material.

Computations in dimension 3 are too tedious to be carried out by hand. They 
have been automated (see \cite{BidegarayFesquet06c}).

\end{document}